\documentclass[11pt,a4paper]{article}

\usepackage{amsmath}
\usepackage{amssymb}
\usepackage{booktabs}

\usepackage[numbers,sort&compress]{natbib}
\usepackage{etoolbox}
\usepackage{xspace}
\usepackage{enumitem}
\usepackage{float}
\usepackage{placeins}
\usepackage{caption}
\usepackage{graphicx}

\usepackage{authblk} 
\usepackage[margin=1in]{geometry} 
\usepackage{hyperref} 

\title{Uncertainty-Aware Crack Growth Forecasting via Conditional Denoising Diffusion Models for Phase-Field Fracture}

\author[1,2]{Jahnavi Krishna Koda\thanks{Email: jkoda@islander.tamucc.edu}}
\author[2]{S. M. Mallikarjunaiah\thanks{Corresponding author. Email: m.muddamallappa@tamucc.edu}}

\affil[1]{Coastal and Marine Systems Science, Texas A\&M University--Corpus Christi, Corpus Christi, TX, USA}
\affil[2]{Department of Mathematics \& Statistics, Texas A\&M University--Corpus Christi, Corpus Christi, TX, USA}

\date{} 

\begin{document}

\maketitle

\begin{abstract}
The accurate prediction of brittle crack initiation, propagation, and complex topological 
evolution remains computationally prohibitive when utilizing traditional high-fidelity 
phase-field finite element methods. To overcome these computational bottlenecks, a 
physics-informed conditional Denoising Diffusion Probabilistic Model (DDPM) is proposed 
for the full-field spatiotemporal forecasting of fracture evolution across diverse loading 
regimes and energy decomposition methods. The generative architecture is conditioned on 
rolling historical damage states and explicitly derived kinematic proxies---phase-field 
velocity and gradient magnitude---ensuring temporal coherence without non-physical 
artifacts. The model is trained and evaluated on a comprehensive dataset of 6{,}000 
phase-field simulations.

The principal contribution is spatially-localized uncertainty quantification without 
modification to the trained model. Ensemble variance concentrates at crack branching 
junctions ($\sigma_\mathrm{max} = 0.222$ at Y-junction bifurcations; zero high-uncertainty 
pixels in four deterministic propagation cases), while the high-$\sigma$ tail identifies 
high-error predictions with 90\% precision---an 18-fold improvement over random 
selection. One-step crack tip localization achieves sub-pixel accuracy (0.12 px mean 
error) across both held-out validation subsets (shear-star and tension-spect), 
confirming cross-regime generalization. In closed-loop autoregressive rollout over 50 
steps, the DDPM maintains Dice = 0.929 $\pm$ 0.010 while a deterministic U-Net collapses 
to Dice = 0.423 under error accumulation, a 2.2$\times$ gap that establishes the 
value of stochastic re-sampling for long-horizon stability. Per-step inference requires 
approximately 3.6 s on an H100 GPU, approximately 28$\times$ faster than the FEM 
reference and 1{,}000$\times$ slower than a deterministic U-Net a cost that buys 
the stochastic diversity enabling uncertainty quantification. This work establishes 
a probabilistically-calibrated generative surrogate for brittle fracture, with direct 
application to real-time structural health monitoring and risk-aware computational 
mechanics.
\end{abstract}


\vspace{1em}

\section{Introduction}
\label{sec:Introduction}
Crack initiation, propagation and complex topology evolution in brittle solids remains one of the major challenges in computational mechanics. Phase-field models of fracture have established themselves as a versatile simulation framework over the last two decades. Based on the variational theory of Francfort and Marigo \cite{FRANCFORT19981319} who transformed Griffith's energy criterion into a free discontinuity minimization problem, phase-field formulations regularize sharp crack surfaces using a continuous scalar damage variable $\phi \in [0,1]$. This regularization eliminates the need for explicit tracking of cracks; meshes do not need to be remeshed; enrichment is also not necessary. Furthermore, complex topological transitions such as nucleation, branching, merging and kinking can be represented within a unified continuum framework \cite{BOURDIN2000797,Bourdin2008Variational,yoon2021quasi,lee2022finite,manohar2025convergent,manohar2025convergence}.

Phase-field fracture was established computationally via thermodynamically consistent variational formulations of Miehe et al. \cite{MieheIJNME2010} who introduced the spectral decomposition of the strain energy to distinguish tensile and compressive crack driving forces, and the companion operator-split algorithm \cite{MIEHE20102765} that enabled staggered solution schemes with robustness. Dynamic fracture extensions were developed by Borden et al. \cite{BORDEN201277}. Ambati et al. \cite{Ambati2014Review} provided a comprehensive review and introduced a hybrid model to reduce computational cost. Wu \cite{Wu2017Unified} later proposed a unified phase-field theory covering damage mechanics and quasi-brittle failure broadening the applicability of the variational framework. Comprehensive reviews of phase-field models for fracture and fatigue are presented by Li et al. \cite{LI2023109419}.

The variety of applications of the phase-field approach has allowed its application to many multiphysical problems. Phase-field models were employed by Chu et al. \cite{Chu2017ThermalShock}  to study dynamic trajectories under thermal shock. Mart\'inez-pa\~neda et al.  \cite{MartinezPaneda2018Hydrogen} coupled hydrogen transport with fracture evolution to model hydrogen-assisted cracking. Hageman and Mart\'inez-pa\~neda \cite{Hageman2023ElectroChemo} extended the framework to electrochemical-mechanical fracture including crack contained electrolytes and chemical reactions. Zhang et al. \cite{Zhang2020Anisotropy} incorporated elastic anisotropy to capture orientation-dependent resistance to fracture. Alternative nonlocal approaches such as peridynamics \cite{Abdoh2022Peridynamics} offer additional representation capabilities for discontinuities in heterogeneous media.

However, despite these advances, a significant computational bottleneck persists across high-fidelity phase-field methods. In addition to fine spatial discretization needed to resolve the regularization length scale $\ell_0$ due to the nonlinear and non-convex nature of the energy minimization required at each loading increment, high computational costs prohibit practical engineering applications. For example, a single two-dimensional simulation on a moderately resolved grid requires hours of wall-clock time on multicore architectures. This cost effectively prevents real-time structural assessment, large-scale parametric studies, Monte Carlo reliability analyses, and uncertainty-quantified design optimization. Therefore, there is a growing need for surrogate models which can reproduce full field outputs of phase-field solvers at a fraction of the computational cost while maintaining physical fidelity.

Phase-field fracture model is a development from the variational re-formulation of Brittle Fracture by Francfort and Marigo \cite{FRANCFORT19981319}, where they proposed minimising the total elastic strain energy and crack surface area for all possible crack sets. 
Bourdin et al. \cite{BOURDIN2000797} were among the first authors to implement numerically this approach using an elliptical regularization of the crack functional which replaced the discrete crack set by a continuous damage field whose size was determined by a length scale parameter $\ell_0$. This regularized formulation was systematically developed in the monograph of Bourdin et al \cite{Bourdin2008Variational}.

Two seminal works published by Miehe and colleagues \cite{MieheIJNME2010,MIEHE20102765} represent a major step towards establishing the general framework of implementation. In the first paper the authors present a formulation which is thermodynamic consistent together with a method to decompose the elastic energy into its spectral components. 
In the second one, the authors propose a numerical solution strategy based on an operator split-staggered scheme that incorporates a local history field ensuring the irreversibility of the cracks. Borden et al. \cite{BORDEN201277} extended phase-field fracture to the fully dynamic case, i.e., taking into account inertia during crack evolution, and demonstrated in three dimensions how cracks can be propagated using adaptative mesh refinement. 
Ambati et al. \cite{Ambati2014Review} presented a review of all existing models for phase-field brittle fractures, as well as a new fast hybrid formulation. 
Finally, Wu \cite{Wu2017Unified} introduced a unified phase-field theory that includes both damage mechanics and quasi-brittle failure through a geometric crack function described by the phase-field profile.

Other studies have focused on some specific aspects of phase-field modelling related to the present work. For example, Ziaei-Rad et al.\cite{ZIAEIRAD2016304} developed methods to extract explicitly crack paths from the diffuse phase-field representation through non-maximum suppression. 
Gerasimov et al. \cite{DBLP:journals/corr/abs-2005-01332} explored stochastic phase-field modelling in order to compute multiple crack paths and their respective probability distributions due to material heterogeneities or numerical perturbations. 
This effect corresponds exactly to what we want to quantify here. 
We generate our phase-field computations using MOOSE \cite{PERMANN2020100430}, a powerful, open-source, massively parallel multi-physics simulation platform.

Machine learning methods are being used to overcome the computational challenges associated with running high fidelity fracture analyses. As such, machine learning has become a widely accepted methodology for developing surrogates of complex fracture simulations. In terms of discrete crack descriptors, Worthington and Chew \cite{Worthington2022CrackPathML}, utilized machine learning algorithms to identify crack paths in random media through treating fracture as a series of directional decision-making processes. Additionally, Karimian et al. \cite{karimian2020neuralnetworkparticlefiltering} integrated both neural networks with particle filters for probabilistic prediction of crack propagation. Furthermore, Perera and Agrawal \cite{perera2022generalizedmachinelearningframework} presented a generalized framework using graph neural networks for modeling brittle cracks utilizing transfer learning for different materials. Although these methodologies have shown efficiency in reducing computational expense, the methodologies described above utilize limited geometric representation and fail to recreate the complete field of spatial-temporal state variables necessary for subsequent mechanical analysis.

Physics informed methodologies have attempted to address this gap. Zhu et al. \cite{ZHU2025117655} proposed an extended Physics Informed Extreme Learning Machine (PIELM) model for simulating Linear Elastic Fracture Mechanics (LEFM), which approximates displacement fields near crack tip regions by enriching the approximate solution space with known singular factors. Manav et al. \cite{MANAV2024117104} investigated the application of the Deep Ritz Method within the Variational Framework of Phase Field Fracture Models; it was determined that physics informed neural networks could simulate crack nucleation, propagation, branching and coalescence. Despite the capability of optimization-based methods to accurately reproduce results for specific problem configurations, the need to resolve for every new problem configuration limits their applicability for parameter study type investigations.

Neural Operator Architectures represent an alternative to previous methodologies by enabling the ability to learn mappings between function spaces. Lu et al. \cite{Lu2021DeepONet} introduced DeepONet, a neural operator based upon the Universal Approximation Theorem for Operators. Similarly, Li et al. \cite{Li2020FourierNO} proposed the Fourier Neural Operator (FNO), which utilizes spectral domain learning and enables resolution-independent generalization. Finally, Kovachki et al. \cite{JMLR:v24:21-1524} provided a comprehensive review of neural operators. Specifically, within fracture mechanics, Goswami et al. \cite{GOSWAMI2022114587} proposed a physics informed variational DeepONet for predicting crack path trajectories in quasi-brittle materials, where the variational governing energy was incorporated into the training loss. More broadly, Goswami et al. \cite{goswami2022physicsinformeddeepneuraloperator} proposed extensions of physics informed deep neural operator networks to a wide range of applications within solid mechanics.
Very recently, Hamdi and Lejeune \cite{HAMDI2026118526} provided a standardized set of 6,000 phase-field simulations representing various loading conditions and initial crack configurations. Further, Hamdi and Lejeune \cite{HAMDI2026118526} also evaluated baseline comparisons of PINNs, FNOs, and U-Nets. Therefore, the present investigation uses this standardized data-set to extend the baseline comparisons to include conditional diffusion models; thus, representing the first evaluation of generative models on this standardized test bed.

A major disadvantage of all current machine learning surrogate methodologies is the use of deterministic predictions. Most conventional machine learning architectures produce one output for each input. This is fundamentally in conflict with the physics of fracture. At structural bifurcation points, there exist multiple mechanically allowable crack paths; therefore, any deterministic prediction will either smear or randomly select from among those allowable paths. The inability to make a deterministic prediction for all possible solutions represents the motivation behind the generative approach adopted in this research.

Denoising Diffusion Probabilistic Models (DDPM) by Ho et al. \cite{Ho2020DDPM} produce synthetic data by learning to reverse a fixed Gaussian noising process through iterative denoising. The theoretical foundations were laid by Song and Ermon \cite{song2020generativemodelingestimatinggradients} who framed generative modeling as score estimation. In addition to improving the efficiency of sampling in diffusion models and expanding the design space of those models, subsequent works include deterministic implicit sampling for faster denoising by Song et al. \cite{song2022denoisingdiffusionimplicitmodels}, a systematic analysis of diffusion architecture design spaces by Karras et al. \cite{karras2022elucidatingdesignspacediffusionbased}, and classifier-free guidance for improved conditional image synthesis by Ho and Salimans \cite{Ho2022ClassifierFreeDG}. For diffusion models in computer vision, Croitoru et al. \cite{Croitoru2023DiffusionSurvey} provide comprehensive reviews.

Conditional diffusion models have also been extended to a variety of other scientific inference tasks. Zhang et al. \cite{Zhang2024cDDPMSeismic} demonstrated that conditional DDPMs could be used to separate and image seismic diffractions and demonstrate the ability to recover very weak physical signal contributions from complex geologic signatures. Jiang et al. \cite{Jiang2024CrackSegDiff} demonstrated the use of diffusion models to segment cracks in images using multi-modality input; Graham et al. \cite{Graham2023DiffusionOOD} demonstrated that diffusion based feature representations could indicate when an input was out-of-distribution; and Nadkarni et al. \cite{Nadkarni2025MolecularDenoising} incorporated physics informed priors into diffusion based pipelines for molecular denoising. However, all of these examples address static inference tasks related to perception, detection, and reconstruction rather than the dynamic prediction of future states of physical systems.

More recently, researchers have begun extending diffusion models to time dependent partial differential equation (PDE) simulations. Lippe et al. \cite{10.5555/3666122.3669068} introduced PDE-Refiner which produces accurate long-horizon roll-outs for fluid dynamics via a diffusion inspired refinement process that preserves the high frequency spectral components that are lost in most current neural PDE solvers. Cachay et al. \cite{cachay2023dyffusiondynamicsinformeddiffusionmodel} proposed DYffusion -- a dynamics-informed diffusion model for spatiotemporal forecasting that combines the denoising process with the temporal evolution of the physical system being simulated. Shu et al. \cite{SHU2023111972} presented a physics informed diffusion model for reconstructing highly accurate flow fields from incomplete or low fidelity measurement data. Kohl et al. \cite{KOHL2026108641} performed a comprehensive benchmark study of autoregressive conditional diffusion models for simulating turbulent flows, showing that diffusion-based methods can provide equal or superior temporal stability compared to deterministic methods, while additionally producing probability distributions consistent with the statistical properties of the underlying physics.

These studies collectively establish that diffusion models are viable surrogates for complex PDE systems and that their probabilistic nature provides a natural mechanism for uncertainty estimation. However, existing applications largely target fluid dynamics. The application of conditional diffusion models to solid mechanics, and specifically to phase-field fracture, where crack irreversibility, topological bifurcation, and sparse damage fields present distinct challenges, remains underexplored. The present work addresses this gap.

Reliability of uncertainty quantification is necessary for the application of surrogate machine learning models in safety-related engineering. There are two major paradigms used for uncertainty quantification in deep learning. Gal and Ghahramani \cite{pmlr-v48-gal16} utilized Monte Carlo dropout during inference to show that it could be considered as a form of approximate Bayesian inference in deep Gaussian processes. They further demonstrated that this method would provide an efficient way to quantify predictive uncertainty. Lakshminarayanan et al. \cite{Lakshminarayanan2017Ensembles} proposed deep ensembles, where each member of the ensemble was created by training independent versions of the network on the same dataset. When aggregating these results, they found that the uncertainties estimated were very accurate when compared against Bayesian uncertainty estimators.

Psaros et al. \cite{PSAROS2023111902}, in addition to developing the use of uncertainty quantification for physics informed neural networks and neural operators, established a set of standards for evaluating uncertainty quantification methods. Furthermore, they evaluated several different methods utilizing the same set of criteria. Yang et al. \cite{Yang_2022}, demonstrated how to perform uncertainty quantification for deep operator networks. In order to scale their method, they employed randomized priors.

Uncertainty quantification is particularly important in fracture mechanics. The phase field formulation of fracture is inherently non-convex. Moreover, in regions of structural bifurcation small variations in either load or material properties can cause a structure to follow one path versus another path. Gerasimov et al. \cite{DBLP:journals/corr/abs-2005-01332}, have shown this stochastic multiplicity within the context of phase fields by predicting multiple crack paths along with the probability associated with each path. However, deterministic machine learning surrogates, i.e., U-Nets, FNOs and DeepONets etc., cannot capture this multiplicity. Instead of producing an ensemble of possible crack paths, each prediction produced represents some type of average of all possible crack paths, potentially leading to a solution that does not physically exist.

A potential solution to this problem lies within the ability of diffusion models to stochastically sample. Specifically, since each reverse diffusion trajectory starts from independent random variables (i.e. normally distributed noise) and converges to a local minimum within the posterior distribution over crack states, an ensemble of independent predictions made for the same input will create a distribution over all possible fracture outcomes. The variance of this distribution will serve as an uncertainty map that concentrates at locations where there is structural instability. This study examines if this theoretical benefit is achieved through practical experimentation involving uncertainty-based prediction for phase field fracture.

This work introduces a conditional denoising diffusion probabilistic model to perform full-field spatio-temporal forecasting of phase-field crack evolutions. Our main contributions were:

\begin{enumerate}

\item \textbf{Physics-aware conditional diffusion architecture.}

A DDPM conditioned on a history of damage fields, including the temporal development of damage, damage velocity and the magnitude of the spatial gradient of damage, for next-step forecasting of the phase-field state. The architecture utilizes a physics-informed U-net backbone \cite{10.1007/978-3-319-24574-4_28}, which contains two separate conditioning pathways. Features representing geometry are concatenated along the channel axis, while scalar features, such as time embedding, load type and crack-density are embedded by an MLP-projection into the bottleneck.

\item \textbf{Cross-regime generalization and crack-specific evaluation protocol.} 

The model is tested on separate hold-out validation sets for shear-star and tensile spectra yielding 4{,}000 predictions in total and it was seen that the model performs well across all four combinations of loading regime and method for decomposing strain energy (unseen during training). In addition to traditional area-based metrics (such as Dice and IoU), a suite of crack region-specific metrics is introduced including Hausdorff distance, modified Hausdorff distance, crack-tip localization error, crack-length error and topological correctness measured using Betti numbers.

\item \textbf{Spatially localized uncertainty quantification at structural bifurcations.} 
The stochastic reverse diffusion process to generate prediction ensembles without model modification is utilized. Ensembles generated from the model exhibit variance concentrated around crack-branching junctions while maintaining uniform low levels of variance throughout regimes where crack propagation is deterministic. High-$\sigma$ tails can be used to identify high-error predictions with a precision of 90\% (i.e., 18${\times}$ better than random), establishing the uncertainty output as a useful tool for engineers to use when identifying high-risk areas. Calibration of the model's uncertainty output is assessed using Spearman Rank Correlation ($\rho = 0.269$) on crack pixels and $N$-scaling coverage analysis (52.2\% at $N=10$ and 60.0\% at $N=50$).

\item \textbf{Autoregressive rollout stability and baseline comparison.} The model was evaluated in closed loop autoregressive mode (with no teacher forcing) over 50 steps. The DDPM maintained a Dice value of 0.929 $\pm$ 0.010. On the other hand, a deterministic U-Net (trained under exactly the same conditions) collapsed to Dice = 0.423. Thus, it has been demonstrated that stochastic re-sampling prevents the error accumulation that degrades fixed point predictors. Error growth exponents are determined for power-law behavior over three different stochastic seeds.

\item \textbf{Computational characterization.} 
Per-step inference was benchmarked at $3{,}637 \pm 30$ ms on an H100 NVL GPU; this corresponds to being about 28$\times$ faster than the FEM reference for each saved step. The computational overhead associated with generating stochastic samples is estimated to be approximately (${\sim}1{,}000\times$) greater than that required by a deterministic U-Net; however, this overhead enables uncertainty quantification.

\end{enumerate}

\FloatBarrier 

\section{Phase-Field Fracture and Benchmark Dataset}
\label{sec:phasefield_dataset}

This section describes the phase-field fracture formulation underlying the training data, specifies the material and numerical parameters, and details the benchmark dataset and evaluation protocol.

\subsection{Governing equations}
\label{subsec:governing_equations}
Phase-field models replace sharp crack discontinuities with a continuous scalar damage variable $\phi \in [0,1]$, where $\phi = 0$, the model represents an undamaged region of material; when $\phi = 1$ it represents a completely damaged region. A variational approach to modeling fracture was first proposed by Francfort and Marigo \cite{FRANCFORT19981319}; Bourdin et al. \cite{BOURDIN2000797,Bourdin2008Variational}  total energy functional of a cracked body $\Omega \subset \mathbb{R}^2$ is expressed as

\begin{equation}
\mathcal{E}(\mathbf{u}, \phi)
=
\int_{\Omega}
g(\phi)\,
\psi_0\!\left(\boldsymbol{\epsilon}(\mathbf{u})\right)
\, dV
+
\int_{\Omega}
G_c \,
\frac{1}{2}
\left[
\frac{1}{\ell_0}\phi^2
+
\ell_0
\left(\nabla \phi \cdot \nabla \phi\right)
\right]
\, dV
-
\mathcal{P}_{\mathrm{ext}},
\label{eq:total_energy}
\end{equation}
where $\mathbf{u}$ is the displacement field, $\boldsymbol{\epsilon}(\mathbf{u}) = \nabla^s \mathbf{u}$ is the small strain tensor, $\psi_0$ is the strain energy density due to elasticity, $G_c$ is the critical energy release rate, and $\ell_0$ is the length-scale controlling the width of the diffuse crack zone \cite{MIEHE20102765}. The degradation function $g(\phi) = (1-\phi)^2$ couples the damage state to the elastic stiffness. As $\ell_0 \rightarrow 0$, the regularized functional $\Gamma$-converges to the original sharp-crack energy of Griffith \cite{Bourdin2008Variational}.

For a linear elastic material subjected to plane-strain conditions, the strain-energy density is given by
\begin{equation}
\psi_0(\boldsymbol{\epsilon})
=
\frac{1}{2}\lambda \, \mathrm{tr}^2(\boldsymbol{\epsilon})
+
\mu \, \mathrm{tr}\left(\boldsymbol{\epsilon}^2\right),
\label{eq:strain_energy_density}
\end{equation}
where $\lambda$ and $\mu$ are the Lam\'e constants derived from Young's modulus $E$ and Poisson's ratio $\nu$.

To prevent crack growth under compressive loading, the strain energy density is decomposed into tensile and compressive components, $\psi_0 = g(\phi)\psi_0^+ + \psi_0^-$, where only the tensile part $\psi_0^+$ is degraded by the phase field. The dataset employed in this work \cite{HAMDI2026118526} contains simulations generated using three distinct energy decomposition methods:

\begin{itemize}
    \item \textbf{Spectral decomposition} \cite{MIEHE20102765}: The strain tensor is broken down by the eigenvalues and eigenvectors of the strain tensor. The degradation function will be applied only to the positive tensile principal strains. It allows for a physically accurate response under compression-dominated states.

    \item \textbf{Volumetric-deviatoric decomposition} \cite{Ambati2014Review}: The strain tensor is separated into volumetric and deviatoric components. Only the positive volumetric component and the full deviatoric component are degrading. This method is less computationally expensive than spectral decomposition while retaining essential tensile/compressive asymmetry.

    \item \textbf{Star-convex decomposition}: A recent modification of the volumetric-deviatoric approach that introduces calibration parameter $\gamma^{\star} = \sigma_e^- / \sigma_e^+$ to independently control compressive and tensile strengths under multiaxial loading. When $\gamma^{\star}=0$, it reduces to the standard volumetric-deviatoric decomposition. More information can be found in \cite{HAMDI2026118526}. 
\end{itemize}

The mechanical equilibrium and phase-field evolution are solved using the hybrid staggered formulation of Ambati et al. \cite{Ambati2014Review}, which utilizes a linear momentum equation in order to reduce computational cost. Crack irreversibility, $\dot{\phi} \geq 0$, is enforced through use of a local history field variable
\begin{equation}
\mathcal{H}^{+}(\mathbf{x},t)
=
\max_{\tau \in [0,t]}
\psi_0^{+}
\left(
\boldsymbol{\epsilon}(\mathbf{x},\tau)
\right),
\label{eq:history_field}
\end{equation}
following Miehe et al. \cite{MIEHE20102765}, the coupled system is solved by alternating minimization. At each load step, iterate between the displacement and phase-field subproblems.

\subsection{Material properties, boundary conditions, and discretization}
\label{subsec:material_boundary_discretization}

All simulations in the benchmark dataset were generated using an open-source finite element framework with the parameters summarized in Table~\ref{tab:material_params}. The domain is a $2~\mathrm{mm} \times 2~\mathrm{mm}$ square discretized with a uniform structured mesh of $800 \times 800$ bilinear quadrilateral Q1 elements, yielding an element size $h = 0.0025~\mathrm{mm}$. This satisfies the standard resolution requirement $h \leq \ell_0/2$ \cite{MIEHE20102765,BORDEN201277} with $\ell_0 = 0.01~\mathrm{mm}$. All simulations are downsampled to a $128 \times 128$ grid for the machine learning pipeline, corresponding to a spatial resolution of approximately $0.0156~\mathrm{mm}$ per pixel.

\begin{table}[H]
\centering
\caption{Material and simulation parameters from the benchmark dataset \cite{HAMDI2026118526}.}
\label{tab:material_params}
\begin{tabular}{ll}
\hline
\textbf{Parameter} & \textbf{Value} \\
\hline
Young's modulus $E$ & $1000~\mathrm{GPa}$ \\
Poisson's ratio $\nu$ & $0.3$ \\
Critical energy release rate $G_c$ & $1~\mathrm{N/mm}$ \\
Regularization length scale $\ell_0$ & $0.01~\mathrm{mm}$ \\
Domain size & $2~\mathrm{mm} \times 2~\mathrm{mm}$ \\
Mesh resolution & $800 \times 800$ Q1 \\
ML grid resolution & $128 \times 128$ \\
Loading increment $\Delta u$ & $1 \times 10^{-6}~\mathrm{mm}$ \\
\hline
\end{tabular}
\end{table}

Two loading modes are applied:

\begin{itemize}
    \item \textbf{Biaxial tension}: uniform vertical displacement is applied to both Top and Bottom Boundaries to drive Mode~I Opening Fracture. There are 5000 increments of loading and snapshots are taken every 50th Step; so that there are 100 Saved States for Each Simulation. Average Wall-Clock Time per Simulation is around 2.8 Hours on a 16-Core Processor.

    \item \textbf{Shear}: Horizontal displacement is applied to the top boundary whilst the bottom boundary remains fixed, thus driving Mode~II sliding fracture. To prevent the formation of fractures that run parallel to the boundaries (and therefore do not provide useful information), dirichlet conditions $\phi=0$ are imposed along the top and bottom edges. The number of increments of shear loading is 10000 and snapshots are taken every 100 steps. As such, as well as providing 100 saved states per simulation; the average wall-clock time required to perform each simulation is also approximately 5.3 hours on a 16-core processor.

\end{itemize}

Initial crack configurations are generated stochastically: each domain contains $n \in [10,20]$ randomly positioned and oriented line cracks of fixed length $l = 0.25~\mathrm{mm}$. Pre-existing cracks are introduced via the initial strain history field formulation of Borden et al. \cite{BORDEN201277}. This random initialization produces substantial diversity in crack morphology, ranging from relatively sparse configurations with few interacting cracks to dense networks exhibiting extensive branching and coalescence.

\subsection{Benchmark dataset}
\label{subsec:benchmark_dataset}

The complete dataset, introduced by Hamdi and Lejeune \cite{HAMDI2026118526}, comprises 6000 phase-field fracture simulations organized into six subsets by the combination of loading mode and energy decomposition method.

\begin{table}[!htbp]
\centering
\caption{Dataset composition. Each subset contains 1000 simulations with distinct random initial crack configurations. Subset names follow the convention \texttt{loading-decomposition}.}
\label{tab:dataset_composition}
\begin{tabular}{lllll}
\hline
\textbf{Subset} & \textbf{Loading} & \textbf{Energy decomposition} & \textbf{Files} & \textbf{Role} \\
\hline
\texttt{tension-spect} & Biaxial tension & Spectral & 1000 & Validation \\
\texttt{tension-vol} & Biaxial tension & Volumetric-deviatoric & 1000 & Train \\
\texttt{tension-star} & Biaxial tension & Star-convex & 1000 & Train \\
\texttt{shear-spect} & Shear & Spectral & 1000 & Train \\
\texttt{shear-vol} & Shear & Volumetric-deviatoric & 1000 & Train \\
\texttt{shear-star} & Shear & Star-convex & 1000 & Validation \\
\hline
\end{tabular}
\end{table}

All simulations are stored as HDF5 files that contain a $[3,101,128,128]$ three-field channel, corresponding to the phase-field damage $\phi$, horizontal displacement $u_x$, and vertical displacement $u_y$, measured at 101 time steps (the first being the unloading increment) along the $128 \times 128$ grid. In addition, load-displacement curves and initial crack-seed geometries are recorded. Only the phase-field channel $\phi$ is used as the prediction target in this work; the displacement channels are auxiliary.

In Figure~\ref{fig:dataset_gallery}, several examples of final crack shapes obtained from all six datasets are shown. The shear-loaded subsets have been found to have significantly more complex morphology than those subjected to tensile loads. Specifically, shear-loaded cracks show significant diagonally oriented propagation and often include many branching points. Tensile-loaded cracks tend to follow deterministic paths that are primarily vertical and have few branch points. Finally, due to the fact that there are three different ways in which the energy can be decomposed into damage and elastic strain energies, the trained surrogate models must capture the underlying physics of fracture, and cannot simply memorize aspects of the numerical formulation.As noted by Hamdi and Lejeune \cite{HAMDI2026118526}, some of the simulations conducted using these different formulations result in similar crack patterns, while others do so with much less similarity. This difference is thought to reflect differences in the physical behavior of fracture under the influence of different constitutive models.

\begin{center}
\includegraphics[
    width=0.72\linewidth,
    height=0.42\textheight,
    keepaspectratio
]{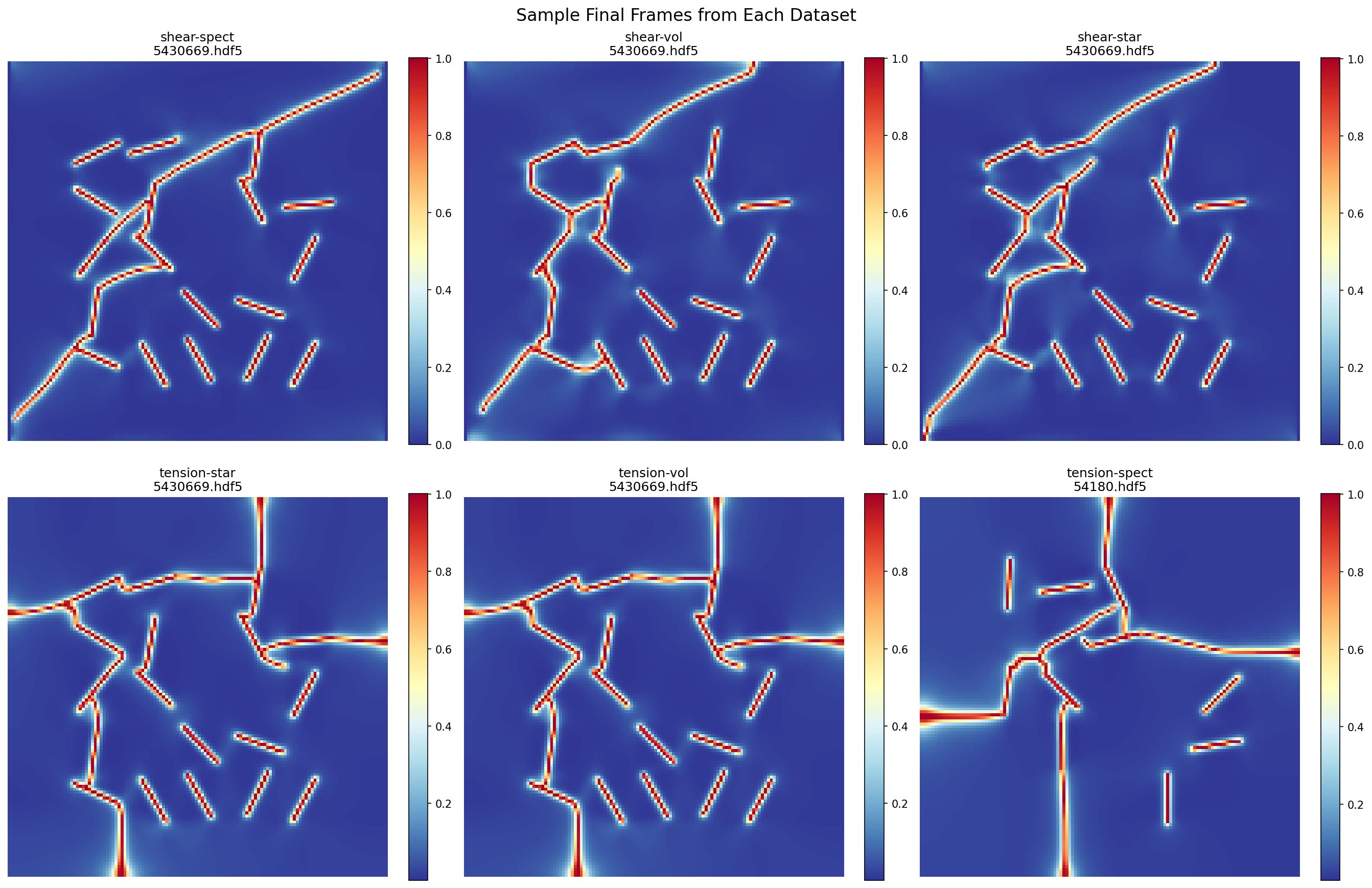}

\captionof{figure}{Representative final-frame crack patterns from the six benchmark subsets. The gallery illustrates the difference between tension- and shear-loaded cases, as well as the influence of spectral, volumetric-deviatoric, and star-convex energy decompositions.}
\label{fig:dataset_gallery}
\end{center}

\subsection{Train/validation split and held-out evaluation protocol}
\label{subsec:train_validation_protocol}

The reason that a subset-based partitioning strategy has been employed here is to test whether the trained surrogate model will generalize to new combinations of loading modes and energy decomposition formulations. There were four training subsets, namely \texttt{shear-spect}, \texttt{shear-vol}, \texttt{tension-star}, and \texttt{tension-vol}. These four training subsets contained a total of 4000 simulations that included both loading modes, and all three of the energy decomposition formulations. On the other hand, the two validation subsets, namely \texttt{shear-star} and \texttt{tension-spect}, had a combined total of 2000 simulations whose specific combinations of loading mode and energy decomposition were never seen before during training.

Therefore, this testing protocol represents a much stricter form of generalization than would be represented by randomly sampling a portion of the simulations at the sample level. That is, the surrogate model is required to extrapolate beyond what was seen during training, specifically with regard to energy decomposition methods. For example, the surrogate model was trained on simulations conducted using shear loading with either spectral or volumetric-deviatoric decompositions of the energy. However, when evaluating its performance on shear loading with the star-convex decomposition method, it is being asked to predict crack patterns that it has never seen before. Therefore, this approach prevents the model from becoming overly specialized in features of the discretization that may not be relevant to other formulations, and instead encourages it to develop physically invariant representations of how cracks evolve.

The surrogate prediction problem is defined such that for each temporal snapshot at time step $t \in [5,100]$, one training sample is generated. At each time step $t$, the surrogate model receives a complete set of conditioning inputs from previous time steps, and produces a predicted value for the phase-field frame at time step $t$. Restricting $t \geq 5$ ensures that there is enough history to generate the necessary input vectors for the conditioning framework outlined in Section~\ref{subsec:physics_aware_features}. Using this definition, a total of 384000 training samples ($4 \times 1000 \times 96$ usable timesteps) are produced, and a total of 192000 validation samples ($2 \times 1000 \times 96$) are produced.

No data augmentations (i.e., no rotation or flipping) are made on the training data. The large number of variations present in terms of both initial crack geometry and loading direction/velocity ensure that there is sufficient variability to avoid overfitting. The excellent match between the training curve and validation curve observed near convergence (as reported in Section~\ref{subsec:training_procedure} supports this conclusion.

\FloatBarrier


\section{Methodology}
\label{sec:methodology}

This section describes the conditional denoising diffusion framework proposed for phase-field crack growth forecasting. The pipeline is organized as follows: Section~\ref{subsec:surrogate_prediction_task} formalizes the surrogate prediction task; Section~\ref{subsec:physics_aware_features} defines the physics-aware conditioning features; Section~\ref{subsec:conditional_ddpm} presents the diffusion framework; Section~\ref{subsec:unet_architecture} details the U-Net architecture; Section~\ref{subsec:training_procedure} describes the training procedure; Sections~\ref{subsec:stochastic_inference} and~\ref{subsec:ensemble_uq} present the stochastic inference and ensemble uncertainty quantification protocols; and Section~\ref{subsec:physical_constraints} addresses the enforcement of physical constraints.

\begin{center}
    \centering
    \includegraphics[
    width=0.82\linewidth,
    height=0.52\textheight,
    keepaspectratio
]{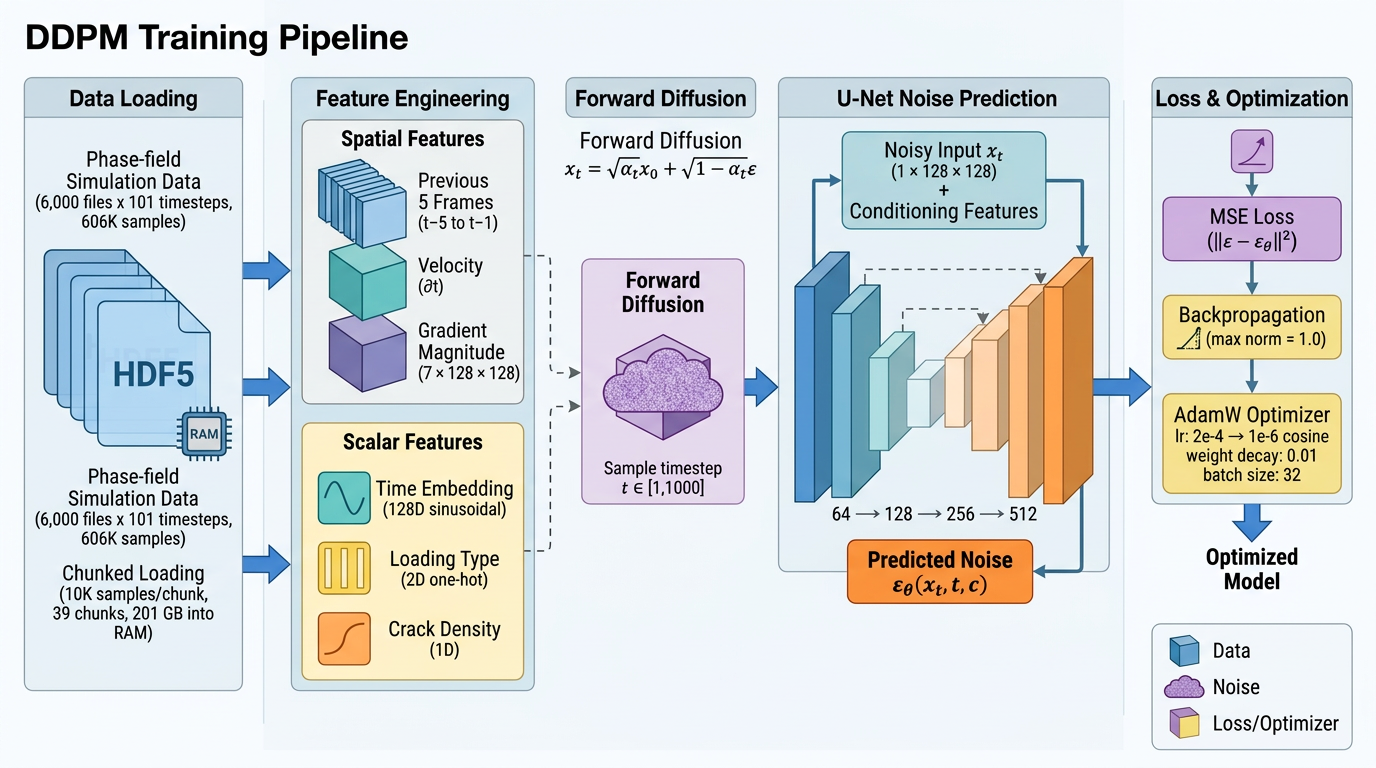}
    \captionof{figure}{Schematic of the conditional U-Net denoising architecture. Spatial conditioning is introduced through channel concatenation, while scalar conditioning is injected through an MLP embedding at the bottleneck.}
    \label{fig:unet_architecture}
\end{center}

\subsection{Surrogate prediction task}
\label{subsec:surrogate_prediction_task}

The objective is to learn a mapping from the recent history of the phase-field damage state to the next temporal frame. Specifically, given a sequence of preceding phase-field snapshots and derived kinematic descriptors at timestep $t$, the model predicts the full-field damage state $\phi(\cdot,t)$ on the $128 \times 128$ spatial grid. The prediction target is a single-channel image $\mathbf{x}_0 \in \mathbb{R}^{1 \times 128 \times 128}$ representing the normalized phase field at timestep $t$, where values are mapped to $[-1,1]$ via the affine transformation
\begin{equation}
\phi_{\mathrm{norm}} = 2\phi - 1.
\label{eq:phi_normalization}
\end{equation}

The model operates in a one-step prediction mode: for each target timestep $t \in [5,100]$, it receives ground-truth conditioning information from timesteps $t-5$ through $t-1$ and predicts the state at $t$. During autoregressive rollout, evaluated in Section~\ref{subsec:rollout_stability}, the model's own predictions replace ground-truth conditioning for multi-step forecasting.

\subsection{Physics-aware feature representation}
\label{subsec:physics_aware_features}

A central design choice is to condition the diffusion model not on raw images alone, but on a feature tensor that explicitly encodes the kinematic drivers of crack propagation. For each target timestep $t$, the conditioning information $\mathcal{F}(t)$ comprises a spatial tensor $\mathbf{F}_{\mathrm{spatial}} \in \mathbb{R}^{7 \times 128 \times 128}$ and a scalar vector $\mathbf{f}_{\mathrm{scalar}} \in \mathbb{R}^{131}$.

\paragraph{Spatial conditioning features.}
Crack propagation is inherently history-dependent: the accumulated stress relaxation field dictates future trajectories. The spatial conditioning tensor contains three categories of information.

First, the five preceding phase-field frames form a temporal history window:
\begin{equation}
\mathbf{H}(t)
=
[
\phi(t-5),\,
\phi(t-4),\,
\phi(t-3),\,
\phi(t-2),\,
\phi(t-1)
],
\label{eq:history_window}
\end{equation}
which captures damage accumulation, growth direction, and local momentum. For $t < 5$, unavailable frames are zero-padded.

Second, the rate of damage evolution is approximated via finite difference:
\begin{equation}
\mathbf{v}(x,y,t)
=
\frac{
\phi(x,y,t-1)-\phi(x,y,t-2)
}{\Delta t},
\label{eq:damage_velocity}
\end{equation}
isolating the active crack-tip zones where the damage front is advancing most rapidly.

Third, the gradient of the phase field at $t-1$ serves as a proxy for local stress intensity:
\begin{equation}
\mathbf{g}(x,y,t)
=
\|\nabla \phi(x,y,t-1)\|_2
=
\sqrt{
\left(\frac{\partial \phi}{\partial x}\right)^2
+
\left(\frac{\partial \phi}{\partial y}\right)^2
},
\label{eq:gradient_magnitude}
\end{equation}
where spatial derivatives are evaluated via central difference schemes. High gradient magnitudes identify the narrow transition zone at crack fronts where the phase field transitions sharply from intact to fractured material.

\paragraph{Scalar conditioning features.}
Global structural state variables are encoded as a 131-dimensional vector.

First, a 128-dimensional sinusoidal positional encoding $\mathbf{e}_{\mathrm{time}}(t)$, following the Transformer convention, smoothly represents the normalized timestep $t/101$.

Second, a one-hot encoding $\mathbf{e}_{\mathrm{loading}} \in \{[1,0]^{\top},[0,1]^{\top}\}$ distinguishes shear from biaxial tension loading, avoiding ordinal bias.

Third, a crack-density scalar
\begin{equation}
\rho_{\mathrm{crack}}(t)
=
\frac{1}{N}
\sum_{x,y}
\mathbb{1}
[
\phi(x,y,t-1) > 0.6
]
\label{eq:crack_density}
\end{equation}
tracks the macroscopic saturation state, signaling whether the domain is in the early nucleation, active propagation, or late saturation regime.

All spatial features are clipped to $[0,1]$ and affine-mapped to $[-1,1]$ to match the diffusion target range. The complete feature schema is summarized in Table~\ref{tab:features}.

\begin{table*}[t]
\centering
\caption{Physics-aware conditioning features used for diffusion model input.}
\label{tab:features}
\begin{tabular}{llllp{4.0cm}}
\hline
\textbf{Feature} & \textbf{Type} & \textbf{Dimension} & \textbf{Range} & \textbf{Physical role} \\
\hline
Previous 5 frames & Spatial & $[5,128,128]$ & $[-1,1]$ & Damage history and momentum \\
Velocity field & Spatial & $[1,128,128]$ & $[-1,1]$ & Growth rate $\partial \phi/\partial t$ \\
Gradient magnitude & Spatial & $[1,128,128]$ & $[-1,1]$ & Stress intensity proxy $\|\nabla \phi\|$ \\
Time embedding & Scalar & $[128]$ & $[-1,1]$ & Temporal position \\
Loading type & Scalar & $[2]$ & $\{0,1\}$ & Shear or tension \\
Crack density & Scalar & $[1]$ & $[0,1]$ & Saturation state \\
\hline
Total spatial input & -- & $7 \times 128 \times 128$ & -- & -- \\
Total scalar input & -- & $131$ & -- & -- \\
Target output & Spatial & $[1,128,128]$ & $[-1,1]$ & Next frame $\phi(t)$ \\
\hline
\end{tabular}
\end{table*}

The inclusion of explicitly derived kinematic features, namely velocity and gradient magnitude, rather than relying on purely end-to-end learning from raw frames, is motivated by the ablation study in Section~\ref{subsec:ablation_studies}, which demonstrates that removing these features causes a 30--40\% increase in prediction error.

\begin{center}
    \centering
    \includegraphics[width=0.85\linewidth]{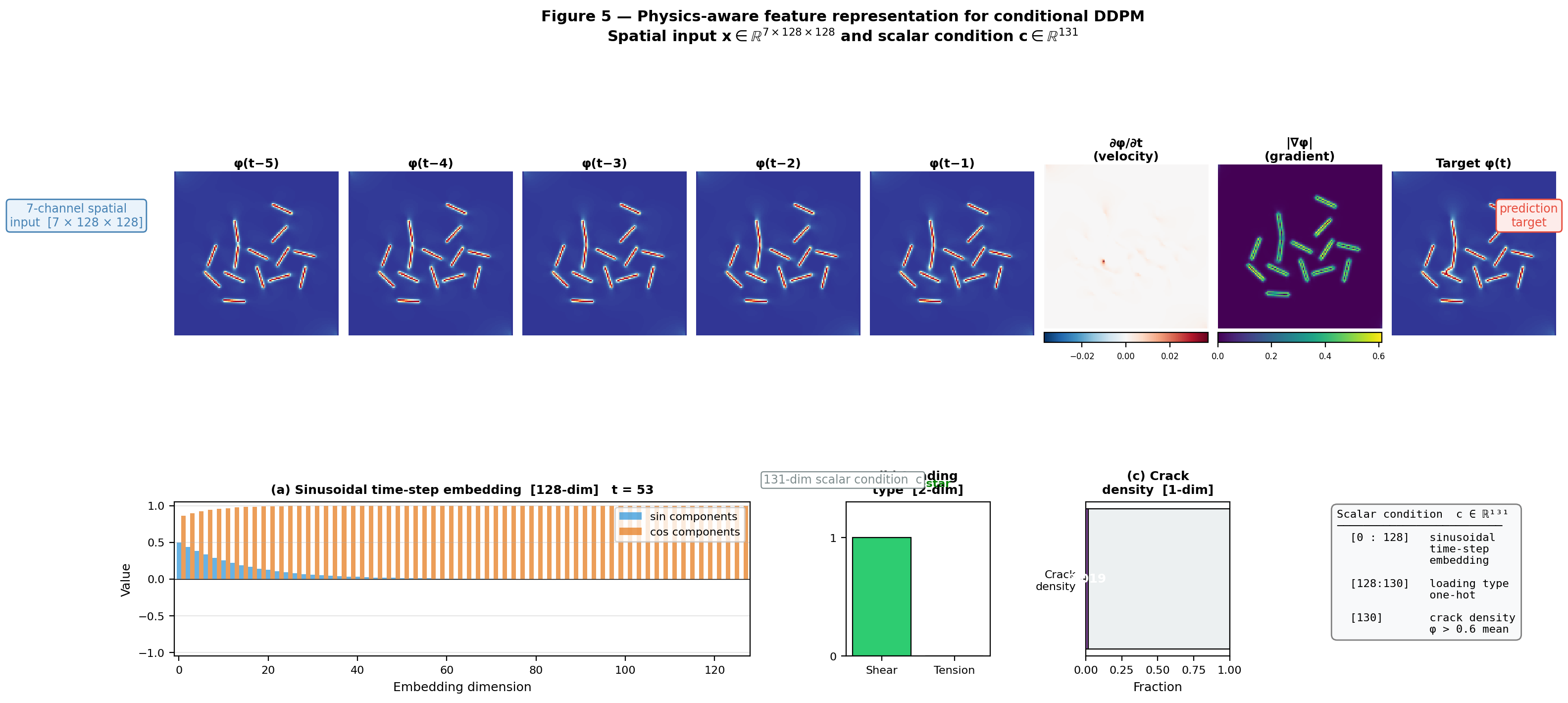}
    \captionof{figure}{Overview of the physics-aware feature representation used to condition the diffusion model.}
    \label{fig:feature_representation}
\end{center}

\subsection{Conditional denoising diffusion framework}
\label{subsec:conditional_ddpm}

The prediction task is formulated as conditional generation using a Denoising Diffusion Probabilistic Model (DDPM) \cite{Ho2020DDPM}. The model synthesizes the target phase-field frame by learning to reverse a fixed Gaussian perturbation process, conditioned on the physics-aware features described in Section~\ref{subsec:physics_aware_features}.

\paragraph{Forward process.}
The forward diffusion progressively corrupts the target frame $\mathbf{x}_0$ over $T=1000$ steps according to a Markov chain:
\begin{equation}
q(\mathbf{x}_k \mid \mathbf{x}_{k-1})
=
\mathcal{N}
\left(
\mathbf{x}_k;
\sqrt{1-\beta_k}\,\mathbf{x}_{k-1},
\beta_k \mathbf{I}
\right),
\label{eq:forward_process}
\end{equation}
where $\beta_k$ follows a linear variance schedule:
\begin{equation}
\beta_k
=
\beta_{\min}
+
(\beta_{\max}-\beta_{\min})
\frac{k}{T},
\qquad
\beta_{\min}=10^{-4},
\qquad
\beta_{\max}=0.02.
\label{eq:beta_schedule}
\end{equation}

Defining $\alpha_k = 1-\beta_k$ and $\bar{\alpha}_k = \prod_{i=1}^{k}\alpha_i$, the noisy state at any arbitrary step $k$ can be sampled in closed form via the reparameterization trick:
\begin{equation}
\mathbf{x}_k
=
\sqrt{\bar{\alpha}_k}\,\mathbf{x}_0
+
\sqrt{1-\bar{\alpha}_k}\,\boldsymbol{\epsilon},
\qquad
\boldsymbol{\epsilon}
\sim
\mathcal{N}(\mathbf{0},\mathbf{I}).
\label{eq:closed_form_noising}
\end{equation}

\paragraph{Reverse process.}
The generative model learns to denoise $\mathbf{x}_k$ back toward $\mathbf{x}_0$, conditioned on the physics-aware context $\mathbf{c}$:
\begin{equation}
p_{\theta}(\mathbf{x}_{k-1} \mid \mathbf{x}_k,\mathbf{c})
=
\mathcal{N}
\left(
\mathbf{x}_{k-1};
\boldsymbol{\mu}_{\theta}(\mathbf{x}_k,k,\mathbf{c}),
\beta_k \mathbf{I}
\right).
\label{eq:reverse_process}
\end{equation}

Following Ho et al. \cite{Ho2020DDPM}, the neural network $\boldsymbol{\epsilon}_{\theta}(\mathbf{x}_k,k,\mathbf{c})$ is trained to predict the noise $\boldsymbol{\epsilon}$ added at step $k$. The posterior mean is then recovered as
\begin{equation}
\boldsymbol{\mu}_{\theta}(\mathbf{x}_k,k,\mathbf{c})
=
\frac{1}{\sqrt{\alpha_k}}
\left[
\mathbf{x}_k
-
\frac{1-\alpha_k}{\sqrt{1-\bar{\alpha}_k}}
\,
\boldsymbol{\epsilon}_{\theta}(\mathbf{x}_k,k,\mathbf{c})
\right].
\label{eq:posterior_mean}
\end{equation}

\paragraph{Training objective.}
The model minimizes the simplified variational bound, which reduces to mean squared error on the predicted noise:
\begin{equation}
\mathcal{L}_{\mathrm{simple}}
=
\mathbb{E}_{k,\mathbf{x}_0,\boldsymbol{\epsilon}}
\left[
\left\|
\boldsymbol{\epsilon}
-
\boldsymbol{\epsilon}_{\theta}
\left(
\sqrt{\bar{\alpha}_k}\,\mathbf{x}_0
+
\sqrt{1-\bar{\alpha}_k}\,\boldsymbol{\epsilon},
k,
\mathbf{c}
\right)
\right\|_2^2
\right],
\label{eq:ddpm_loss}
\end{equation}
where $k$ is sampled uniformly from $\{1,\ldots,T\}$ during training.

\subsection{U-Net architecture and conditioning mechanism}
\label{subsec:unet_architecture}

The noise prediction network $\boldsymbol{\epsilon}_{\theta}$ is parameterized by a U-Net \cite{10.1007/978-3-319-24574-4_28} with a base channel width of 64 and a total of 36,266,113 trainable parameters. The architecture follows the standard encoder-decoder structure with skip connections, adapted for conditional generation through two parallel conditioning pathways.

\paragraph{Spatial conditioning through channel concatenation.}
The 7-channel spatial feature tensor is concatenated with the noisy input $\mathbf{x}_k$ along the channel axis, producing an 8-channel input:
\begin{equation}
\mathbf{x}_{\mathrm{cond}}
=
\mathrm{concat}
\left(
[
\mathbf{x}_k,\,
\phi_1,\,
\phi_2,\,
\phi_3,\,
\phi_4,\,
\phi_5,\,
\mathbf{v},\,
\mathbf{g}
]
\right)
\in
\mathbb{R}^{8 \times 128 \times 128}.
\label{eq:spatial_concat}
\end{equation}

This concatenation preserves exact spatial correspondence between the conditioning fields and the prediction target, ensuring that local crack-tip features in the velocity and gradient fields are spatially aligned with the regions they condition.

\paragraph{Scalar conditioning through MLP injection.}
The 131-dimensional scalar vector is mapped through a multilayer perceptron into a 512-dimensional dense embedding:
\begin{equation}
\mathbf{c}_{\mathrm{embed}}
=
\mathrm{MLP}
\left(
[
\mathbf{e}_{\mathrm{time}},\,
\mathbf{e}_{\mathrm{loading}},\,
\rho_{\mathrm{crack}}
]
\right)
\in
\mathbb{R}^{512}.
\label{eq:scalar_embedding}
\end{equation}

This embedding is injected into the U-Net bottleneck layers via broadcast addition, modulating the feature maps globally based on the macroscopic damage state and loading configuration. The diffusion timestep $k$ is encoded separately using a standard sinusoidal embedding and added to the scalar conditioning pathway, following the convention established in Ho et al. \cite{Ho2020DDPM}.

\subsection{Training procedure}
\label{subsec:training_procedure}

\paragraph{Data loading strategy.}
The full training corpus, consisting of 384,000 samples from four dataset subsets, totals approximately 201~GB in raw form, exceeding typical GPU-host memory budgets for concurrent loading. We address this via a chunked preprocessing pipeline: pre-extracted feature tensors are serialized into $C=39$ chunks of 10,000 samples each, with each chunk occupying approximately 5.2~GB. During training initialization on high-memory compute nodes with 256~GB RAM, all chunks are loaded into volatile memory, bypassing disk I/O entirely. This increases training throughput from 0.4 to 10.9 iterations/s, a factor of $27\times$.

\paragraph{Optimization.}
The model is trained with the AdamW optimizer \cite{Loshchilov2017DecoupledWD} with global gradient clipping, $\|\nabla_{\theta}\|_2 \leq 1.0$, to stabilize early training when the noise prediction residuals are large. The learning rate follows a cosine annealing schedule without restarts:
\begin{equation}
\eta(e)
=
\eta_{\min}
+
\frac{1}{2}
(\eta_{\max}-\eta_{\min})
\left[
1+
\cos
\left(
\frac{\pi e}{E}
\right)
\right],
\label{eq:cosine_lr}
\end{equation}
decaying from $\eta_{\max}=2 \times 10^{-4}$ to $\eta_{\min}=10^{-6}$ over $E=50$ epochs with batch size $B=32$. Training uses FP16 mixed precision via PyTorch's automatic mixed-precision module on a single NVIDIA H100 NVL GPU with 95~GB VRAM. The complete training configuration is summarized in Table~\ref{tab:training}.

\begin{table}[t]
\centering
\caption{Training hyperparameters.}
\label{tab:training}
\begin{tabular}{ll}
\hline
\textbf{Hyperparameter} & \textbf{Value} \\
\hline
Optimizer & AdamW \\
Initial learning rate $\eta_{\max}$ & $2 \times 10^{-4}$ \\
Final learning rate $\eta_{\min}$ & $10^{-6}$ \\
Learning-rate schedule & Cosine annealing, $T_{\max}=100$ \\
Batch size & 32 \\
Epochs & 50 \\
Gradient clipping & Max norm $=1.0$ \\
Mixed precision & FP16 \\
Diffusion steps $T$ & 1000 \\
Variance schedule & Linear, $\beta_{\min}=10^{-4}$, $\beta_{\max}=0.02$ \\
Hardware & NVIDIA H100 NVL, 95~GB \\
\hline
\end{tabular}
\end{table}

\paragraph{Convergence.}
As shown in Figure~\ref{fig:training_convergence}, training achieves stable convergence by approximately epoch 30. The final noise-prediction losses are $\mathcal{L}_{\mathrm{train}} = 6.72 \times 10^{-5}$ and $\mathcal{L}_{\mathrm{val}} = 6.88 \times 10^{-5}$ at epoch 49, differing by less than 2.5\%. This close agreement indicates that the combination of physics-aware conditioning and the dataset-level train/validation split does not induce overfitting despite the large number of model parameters, 36.3M, relative to the task complexity. Total training time is approximately 15.3 GPU-hours.

We note that these loss values represent the noise-prediction objective in Eq.~\eqref{eq:ddpm_loss} and should not be compared directly with the final phase-field reconstruction MSE reported in Section~\ref{subsec:one_step_accuracy}, which measures the error of the fully denoised output.

\begin{center}
    \centering
    \includegraphics[width=0.85\linewidth]{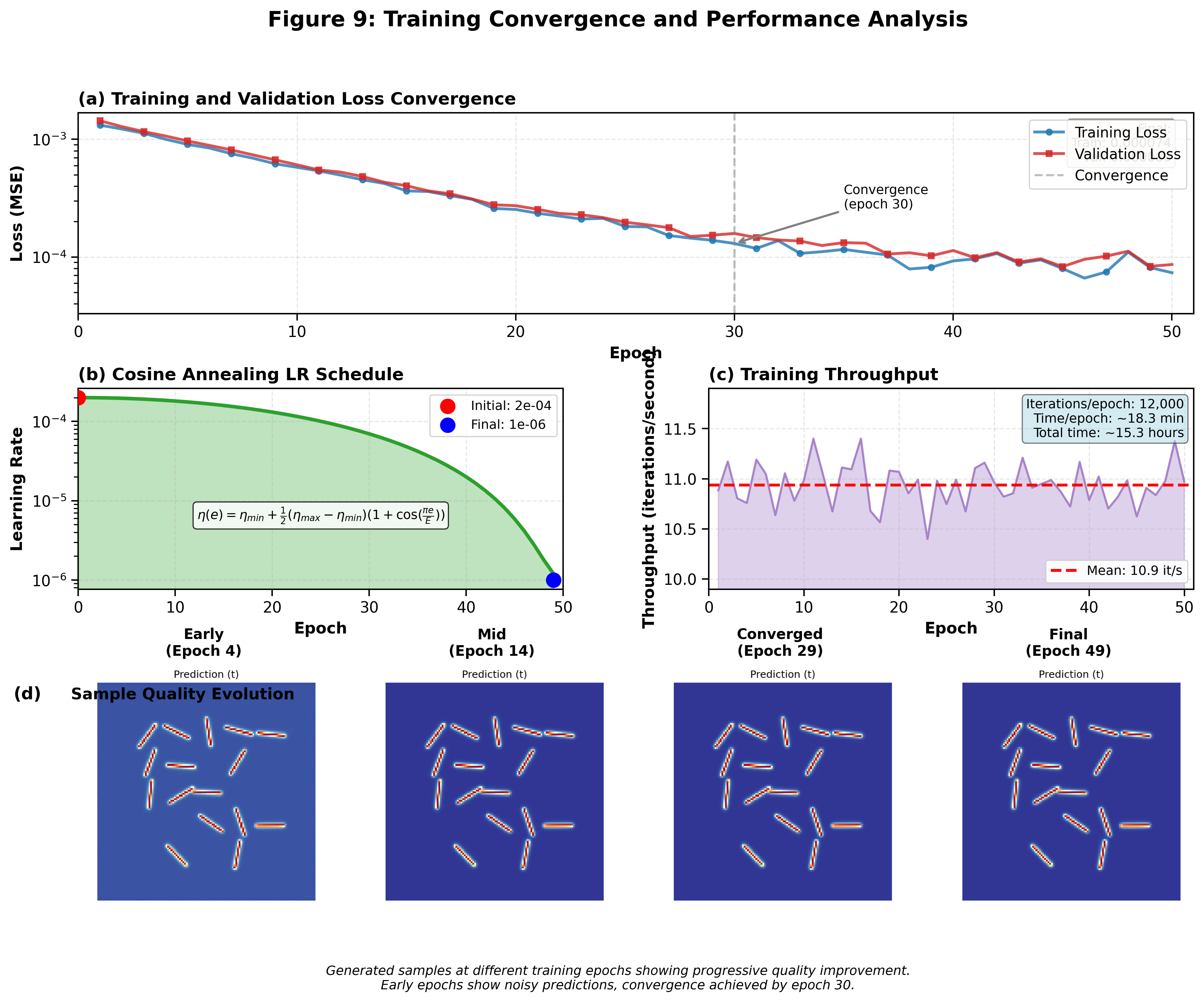}
    \captionof{figure}{Training and validation noise-prediction loss curves for the conditional DDPM.}
    \label{fig:training_convergence}
\end{center}

\subsection{Stochastic inference}
\label{subsec:stochastic_inference}

At inference, the model generates a prediction by executing the complete reverse diffusion chain from standard Gaussian noise $\mathbf{x}_T \sim \mathcal{N}(\mathbf{0},\mathbf{I})$:
\begin{equation}
\mathbf{x}_{k-1}
=
\frac{1}{\sqrt{\alpha_k}}
\left[
\mathbf{x}_k
-
\frac{1-\alpha_k}{\sqrt{1-\bar{\alpha}_k}}
\,
\boldsymbol{\epsilon}_{\theta}(\mathbf{x}_k,k,\mathbf{c})
\right]
+
\sqrt{\beta_k}\,\mathbf{z},
\qquad
\mathbf{z}
\sim
\mathcal{N}(\mathbf{0},\mathbf{I}),
\label{eq:reverse_sampling}
\end{equation}
for $k=T,T-1,\ldots,1$. The stochastic noise injection $\mathbf{z}$ is retained at each step; that is, we do not use the deterministic DDIM sampler \cite{song2022denoisingdiffusionimplicitmodels}, as this stochasticity is essential for generating diverse samples that capture the multimodal posterior over crack paths.

The full 1000-step reverse chain is required. We verified experimentally that reducing to 100 steps produces catastrophically degraded predictions, with the Dice coefficient dropping from 0.995 to 0.160. This indicates that the model's learned denoising trajectory cannot be aggressively truncated without substantial retraining or architectural changes such as progressive distillation.

Single-frame inference requires approximately 3.6\,s on an NVIDIA H100 NVL GPU 
(3{,}637\,ms on a cold GPU at 29$^\circ$C; timing varies by ${\sim}15\%$ with 
thermal state), with a peak memory footprint of 488\,MB. The output is rescaled from $[-1,1]$ back to $[0,1]$ via
\begin{equation}
\phi = \frac{\mathbf{x}_0 + 1}{2}
\label{eq:phi_rescale}
\end{equation}
before evaluation.

\subsection{Ensemble uncertainty quantification}
\label{subsec:ensemble_uq}

The stochastic nature of the reverse diffusion process provides a natural mechanism for uncertainty quantification without any modification to the trained model. Each independent execution of the reverse chain from a different initial noise sample $\mathbf{x}_T^{(i)}$ produces a distinct but physically plausible realization of the predicted damage field.

For a given structural condition, namely a fixed set of conditioning features $\mathbf{c}$, we generate an ensemble of $N$ independent predictions $\{\mathbf{x}_0^{(1)},\ldots,\mathbf{x}_0^{(N)}\}$ and compute the spatially resolved ensemble mean and standard deviation:
\begin{equation}
\bar{\mathbf{x}}_0(x,y)
=
\frac{1}{N}
\sum_{i=1}^{N}
\mathbf{x}_0^{(i)}(x,y),
\qquad
\sigma(x,y)
=
\sqrt{
\frac{1}{N-1}
\sum_{i=1}^{N}
\left(
\mathbf{x}_0^{(i)}(x,y)
-
\bar{\mathbf{x}}_0(x,y)
\right)^2
}.
\label{eq:ensemble_statistics}
\end{equation}

The standard deviation field $\sigma(x,y)$ serves as a pixel-wise uncertainty map. In a well-calibrated model, this map should exhibit low values in regions where the crack path is mechanically determined, such as straight propagation under uniaxial tension, and elevated values at locations of structural instability where multiple crack paths are physically admissible, such as branching junctions.

This approach is conceptually analogous to the deep ensemble method of Lakshminarayanan et al. \cite{Lakshminarayanan2017Ensembles}, but with a key distinction: whereas deep ensembles require training multiple independently initialized models, diffusion-based ensembling exploits the inherent stochasticity of a single trained model. The ensemble size is set to $N=10$ for the calibration analysis in Section~\ref{subsec:calibration_validation} and $N=3$ for the autoregressive rollout in Section~\ref{subsec:rollout_stability}, balancing calibration quality against computational cost. At approximately 3{,}550 ms per sample, a 10-sample ensemble completes in approximately 
35.5 s, remaining substantially faster than the FEM reference time of 2--5 hours per simulation.

Ensemble convergence was verified by comparing $\sigma$ maps at $N=5$, $10$, and $20$ for representative stochastic cases; the uncertainty hotspot location and magnitude stabilize beyond $N=10$.

\begin{center}
    \centering
    \includegraphics[width=0.85\linewidth]{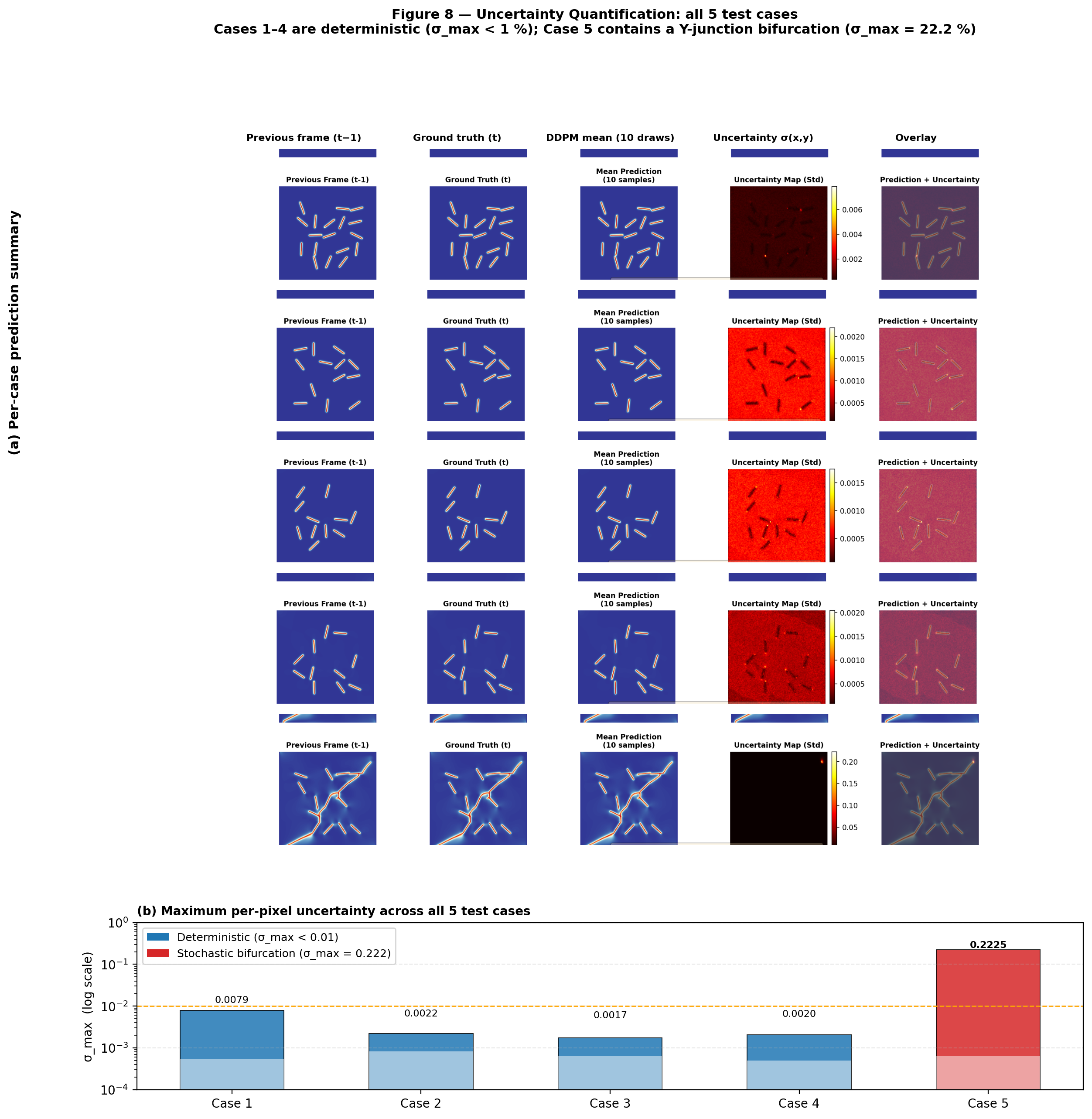}
    \captionof{figure}{Diffusion-based ensemble uncertainty quantification. Multiple stochastic reverse trajectories produce distinct crack-field samples, from which the ensemble mean and pixel-wise uncertainty map are computed.}
    \label{fig:ensemble_uq}
\end{center}

\subsection{Physical constraint enforcement}
\label{subsec:physical_constraints}

Phase-field fracture solutions must satisfy two physical constraints.

\paragraph{Boundedness.}
The damage variable is defined on $\phi \in [0,1]$. The diffusion model operates in the rescaled space $[-1,1]$, and the output is clipped to this range before rescaling. In practice, boundedness is learned implicitly by the model: across the full evaluation set of 2000 predictions, zero out-of-range values were observed, as discussed in the earlier section. 

\paragraph{Irreversibility.}
In quasi-static brittle fracture, damage is monotonically non-decreasing: $\dot{\phi} \geq 0$. This reflects the thermodynamic requirement that cracks cannot heal. While the diffusion model does not enforce this constraint architecturally, we apply a post-hoc monotonic projection:
\begin{equation}
\phi_{\mathrm{corrected}}(x,y,t)
=
\max
\left(
\phi_{\mathrm{pred}}(x,y,t),
\phi(x,y,t-1)
\right),
\label{eq:monotonic_projection}
\end{equation}
which replaces any predicted value that is lower than the preceding frame with the previous value. This projection is physically exact because it enforces the constraint without introducing bias and is computationally negligible. As reported in the earlier section, irreversibility violations are rare, affecting 0.16\% of pixels, and are of small magnitude, approximately $0.0015$. The monotonic projection simultaneously enforces the constraint and reduces the mean absolute error by approximately 6\%.

The choice to enforce irreversibility via post-processing rather than architectural constraint, such as restricting the output activation, is deliberate: a hard architectural constraint would alter the diffusion model's training dynamics and potentially interfere with the stochastic sampling that enables uncertainty quantification. The post-hoc approach preserves the generative model's expressiveness during training while guaranteeing physical compliance at inference.

\FloatBarrier 


\section{Evaluation Framework}
\label{sec:evaluation}
This area describes how metrics, baselines, and procedures were developed to assess the conditional DDPM. Evaluation will be conducted across 5 categories: crack-specific precision (Section~\ref{subsec:crack_metrics}), comparison with established surrogate architectures (Section~\ref{subsec:baseline_models}), long-horizon autoregressive stability (Section~\ref{subsec:autoregressive_rollout_protocol}), uncertainty calibration (Section~\ref{subsec:uncertainty_calibration}), and computational efficiency (Section~\ref{subsec:computational_benchmarking}).

\subsection{Crack-specific topological metrics}
\label{subsec:crack_metrics}
Pixel-level global metrics such as mean squared error (MSE), and pixel accuracy are greatly influenced by the undamaged background in phase-field fracture simulations. In these types of simulations, cracks usually account for less than 2-4\% of the total spatial domain. Hamdi and Lejeune \cite{HAMDI2026118526} also stated that due to the large amount of imbalance in the classes present in the data, pixel accuracy was rendered un-informative. For example, a model that predicts that there is no damage anywhere would have an accuracy of better than 96\%. In fact, on average, cracked pixels represent approximately 2.9\% of the spatial domain across the entire evaluation set. Therefore, due to the dominance of the undamaged background, metrics like Dice are somewhat buoyed up by it; and the more discriminatory assessment comes from boundary-sensitive metrics (e.g. Hausdorff Distance, Tip Location Error). Therefore, the quality of the crack-region specific predictions is assessed using a suite of eight crack-region specific metrics calculated over the binarized crack mask, $\phi = 0$, requiring at least twenty cracked pixels for a valid evaluation.

\paragraph{Overlap metrics.}
The Dice coefficient and intersection-over-union (IoU) calculate the degree to which predicted and ground truth crack masks share common spatial locations:
\begin{equation}
\mathrm{Dice}
=
\frac{2|P \cap G|}{|P| + |G|},
\qquad
\mathrm{IoU}
=
\frac{|P \cap G|}{|P \cup G|},
\label{eq:dice_iou}
\end{equation}
where $P$ and $G$ refer to the sets of pixels that belong to the predicted and ground-truth crack masks, respectively. These metrics both vary between 0 (no overlap) and 1 (perfect agreement), but they do so independently of the overwhelming majority of true negatives found in the undamaged background class.

\paragraph{Geometric distance metrics.}
In order to quantify misalignments beyond simple overlap, two geometric distance metrics are employed based on skeletal representations of the crack masks. The binary crack masks are first skeletonized using \texttt{skimage.morphology.skeletonize}, and then distances are calculated over the resulting one-pixel wide crack curves:

\begin{equation}
\mathrm{Hausdorff}
=
\max
\left(
\max_{p \in S_P}
\min_{g \in S_G}
\|p-g\|,
\;
\max_{g \in S_G}
\min_{p \in S_P}
\|g-p\|
\right),
\label{eq:hausdorff}
\end{equation}
Where $S_P$ and $S_G$ are the sets of points defining the predicted and ground truth skeletons. The Hausdorff distance measures the worst-case misalignment and is reported in pixels. The Modified Hausdorff distance replaces the max with a mean and provides a measure of the average proximity between points in the predicted and ground truth skeletons that is less affected by outliers.

\paragraph{Crack tip location error.}
Crack tips are determined as end-points of the skeleton. End-points are those pixels within the skeleton that have exactly one neighboring pixel within an 8-connectivity neighborhood. The tip location error is the symmetric mean nearest-neighbor or chamfer-distance between the predicted and ground-truth tip-point sets:
\begin{equation}
e_{\mathrm{tip}}
=
\frac{1}{2}
\left(
\frac{1}{|T_P|}
\sum_{p \in T_P}
\min_{g \in T_G}
\|p-g\|
+
\frac{1}{|T_G|}
\sum_{g \in T_G}
\min_{p \in T_P}
\|g-p\|
\right),
\label{eq:tip_error}
\end{equation}
Where $T_P$ and $T_G$ define the sets of points representing predicted and ground-truth tip points. The metric is reported in pixels and represents the capability of the model to correctly locate some of the most mechanically significant features of the damaged field.

\paragraph{Crack length error.}
The absolute difference in total skeleton length between the predicted and ground-truth skeletons provides a measure of the models ability to simulate accurate amounts of damage growth.

\paragraph{Topological correctness through Betti numbers.}
A description of the topological structure of the crack-network can be obtained by computing Betti-numbers derived from the binarized crack-mask. $\beta_0$ counts the number of connected crack components, while $\beta_1$ counts the number of enclosed loops. The latter is computed as
\begin{equation}
\beta_1 = \beta_0 - \chi,
\label{eq:betti_one}
\end{equation}
where $\chi$ is the Euler characteristic. The Betti-0 match rate reports the fraction of test samples where the predicted and ground-truth component counts agree exactly. The Betti-1 match rate is reported conditionally: since closed crack loops, $\beta_1 > 0$, occur in only 0.6\% of test frames, the unconditional match rate is dominated by trivially correct $\beta_1 = 0$ predictions and is not used as a primary metric.

Conditional Betti-1 match rates are reported separately for each validation subset: 
shear-star (12/12 closed loops, 100\%; Wilson 95\% CI: [76\%, 100\%]) and tension-spect 
(60/63 closed loops, 95.2\%; the three failures occur in late-stage frames with dense 
branching topology). 

All metrics are evaluated on 2000 predictions. Each prediction corresponds to 200 withheld simulations with 10 random time-steps selected from each simulation, utilizing 1000 denoising iterations for each prediction with DDPM.

\subsection{Baseline models}
\label{subsec:baseline_models}
Hamdi and Lejeune \cite{HAMDI2026118526} have recently conducted an exhaustive comparative study of several well-known surrogates used to approximate PDE-based simulation solutions for material failure. To provide context for the performance of the conditional DDPM, we will use these previous studies as comparisons. As such, they represent three fundamentally different ways to estimate PDE-based solutions.

\paragraph{Deep Ritz Method.}
Manav et al. \cite{MANAV2024117104} developed a physics-informed approach based on minimizing the variational energy functional of the phase-field problem, using a deep feedforward network that includes 51,109 trainable parameters, six hidden layers with 100 units each and adaptive ReLU activation. In addition, they utilized transfer learning: the weights obtained at one load increment serve as initial weights for subsequent increments. Hamdi and Lejeune \cite{HAMDI2026118526} assessed this approach using their previously-developed dataset and found that no previously-trained models had converged to the correct crack configuration for the complex multiple-crack configurations represented within their dataset. Instead, each randomly-generated seed produced a physically-wrong crack trajectory. Therefore, the Deep Ritz Method represents a "negative" baseline for demonstrating the limitations of contemporary physics-informed methods when applied to complex fracture configurations.

\paragraph{Fourier Neural Operator.}
Li et al. \cite{Li2020FourierNO} proposed a new methodology for modeling functions in the frequency-domain the Fourier Neural Operator (FNO). They accomplished this through parameterization of the integral operator in the frequency-domain. Hamdi and Lejeune \cite{HAMDI2026118526}  trained FNOs having 465,557 trainable parameters, 4 Fourier layers with 12 modes per direction, and GELU activation on all three field channels ($\phi$, $u_x$, and $u_y$), utilizing an unrolled auto-regressive training algorithm with a 10-step history window. Each of the five independently-seeded FNO models used a combination of ground truth values for the first ten time-steps followed by 90 additional time-steps generated autoregressively. Over the five independent models tested, the average Dice score for the test set was 0.621, with an average Dice score of 0.733 achieved via the stacking ensembling technique.

\paragraph{U-Net.}
Ronneberger et al. \cite{10.1007/978-3-319-24574-4_28} presented the U-net encoder-decoder structure, which has been widely adopted as a standard baseline architecture for image-to-image mapping tasks. Hamdi and Lejeune \cite{HAMDI2026118526} trained U-nets with 1,926,689 trainable parameters and ReLU activation on just the phase-field channel using a one-shot non-auto-regressive training paradigm with Dice and focal loss. On average, across five independent seeds, the average Dice score was 0.680 with an average Dice score of 0.733 via stacking ensembling.

\paragraph{Comparison protocol.}
There are two primary forms of comparisons relevant here. When evaluating the predictive accuracy of the DDPM as an auto-regressive predictor (i.e., predicting subsequent values given prior values), it should be compared to the U-Net baseline because both are designed for one-step predictions conditioned upon prior frame information. The Deep Ritz Method baseline provides a reference for the current state of physics-informed approaches. We report the baseline results as published by Hamdi and Lejeune \cite{HAMDI2026118526} on the same dataset; the conditional DDPM is evaluated on the same held-out subsets, \texttt{shear-star} and \texttt{tension-spect}, to ensure a consistent comparison.

There are two caveats regarding the comparison. First, Hamdi and Lejeune's existing studies were trained on only a portion of their benchmark dataset (specifically biaxial tension with spectral decomposition), whereas our conditional DDPM is trained on four portions of their benchmark dataset concurrently. While providing more training data than Hamdi and Lejeune's existing studies, it does require greater generality across loading modes and energy decompositions. Second, there are some minor differences between how Hamdi and Lejeune have evaluated their baselines versus how we evaluate our conditional DDPM. Specifically, Hamdi and Lejeune evaluated their FNOs over a 90-step roll-out from a 10-step seed, whereas we evaluate our conditional DDPM via a 50-step roll-out from a 5-step seed. These differences are noted where relevant in the results discussion.
\subsection{Autoregressive rollout protocol}
\label{subsec:autoregressive_rollout_protocol}
To examine how well long-term predictive performance remains stable, the model is run in closed-loop autoregressive mode (i.e., the model’s output of $\phi$ is substituted into the input window for the next time step $k+1$). Following the first five ground truth frames, no teacher forcing was applied. The process for rolling out the model is shown below:

\begin{enumerate}
    \item \textbf{Seed:} 5 ground-truth frames from a previously unseen simulation.
    \item \textbf{Rollout length:} 50 additional steps past the last frame of the seed.
    \item \textbf{Ensemble per step:} $N=3$ independent DDPM samples were generated. The average of these $N=3$ samples was taken as both the rollout prediction and as the input conditioning for the next time step.
    \item \textbf{Trajectories:} 20 unique simulations were performed and rolled out separately.
    \item \textbf{DDPM steps:} Each prediction took 1000 iterations of DDPM.
\end{enumerate}

Error growth is characterized by fitting the step-dependent error to a power law:
\begin{equation}
\mathrm{error}(k) = a k^b,
\label{eq:power_law_error}
\end{equation}
By performing an ordinary least-squares fit of the data in log-log space using all 50 time steps, it is possible to determine if error growth is sub-linear ($b < 1$), linear ($b = 1$), or super-linear ($b > 1$). In order to estimate the reliability of the fitted exponent, the entire rollout process was repeated three times with different random initialization vectors for the stochastic processes involved. Therefore, b will be presented as the mean $\pm$ std-deviation across trials.

The choice of $N=3$ was made to find a good balance between the quality of rollouts and computational costs. A full rollout of twenty trajectories takes approximately 3.5 hours of GPU time at $N=3$. At $N=10$, this would take about 11.5 hours. Thus, while the main claims regarding predictive stability are based on the absolute Dice score at $k=50$, which is very robust across trials, the estimates of the power-law exponents have larger trial-to-trial variability.

\subsection{Uncertainty calibration assessment}
\label{subsec:uncertainty_calibration}

The analysis of uncertainty quantification measures how accurately the model’s predicted uncertainty reflects its true prediction error. This analysis is composed of four parts.

\paragraph{Pixel-wise calibration.}
At each location on each sample of the test set, $N=10$ ensemble members were generated. From these ten ensembles, the pixel-wise standard deviations $\sigma(x,y)$ were calculated. The regression of the model's predicted uncertainty $\sigma$ versus the absolute difference $|\hat{y}(x,y)-y(x,y)|$ between each ensemble member's prediction and the corresponding ground-truth value, for all x, y and test samples, is then evaluated. Ideally, there will be a strong positive correlation (i.e., $R^2 \rightarrow 1$) between the two quantities and the slope of the regression line will be close to one.

\paragraph{Probabilistic coverage.}
The probability that the ground-truth value lies within a given range (the 95\% prediction interval) around each pixel's predicted values is estimated. The ideal amount of probabilistic coverage is 95\%; amounts significantly above or below this indicate miscalibration.

\paragraph{Reliability diagram.}
Each sample of the test set is divided into deciles according to their predicted uncertainty. The average magnitude of prediction error in each decile is then plotted. An ideally calibrated model should result in a plot whose curves increase monotonically along with the identity line.

\paragraph{Engineering triage precision.}
In addition to statistical calibration, we evaluate how useful the uncertainty metric is for targeted model verification.  If an engineer needs to identify the 10\% of predictions most likely to contain significant errors, defined as any pixel with $|\mathrm{error}| > 0.01$ and compare specifically through selecting either at random, by highest mean uncertainty, or by highest maximum uncertainty among those uncertainties estimated for all samples in a dataset. Precision i.e., what fraction of such selected samples actually contained high-error regions of these verification strategies.

All calibration analyses use 200 unique test cases. These are separate from the 2000 one-step evaluation samples described in Section~\ref{subsec:crack_metrics}, since the uncertainty quantification analysis requires generating ten independent ensemble predictions per sample for a total of 2000 DDPM inference calls.

\subsection{Computational benchmarking protocol}
\label{subsec:computational_benchmarking}

Latency due to Inference is evaluated on an individual NVIDIA H100 NVL GPU equipped with 95~GB of video random access memory (vram), utilizing pytorch version 2.1.0 and cuda version 12.1. A given measurement of latency will consist of 10 iterations used to "warm-up" the gpu prior to the collection of data; these iterations will be removed from the statistical analysis. Following the 10 "warm-up" iterations, there will be 20 iterations collected for use in determining the mean and standard deviation of each sample. Latency will be expressed as mean $\pm$ standard deviation.

Three different inference protocols have been selected for evaluation:

\begin{enumerate}
    \item \textbf{DDPM-1000, proposed:} full 1000-step reverse diffusion chain.
    \item \textbf{DDPM-100, reduced steps:} this configuration utilizes a 100-step reverse diffusion chain. It was selected to assess the costs associated with reducing the number of steps within the reverse diffusion chain, while also assessing whether or not such a strategy could provide sufficient acceleration to make the process useful for the currently proposed model.
    \item \textbf{Deterministic U-Net, single forward pass:} this is a deterministic U-Net that employs the same backbone architecture as the diffusion process being evaluated. Its purpose is to measure the inherent cost of the diffusion process relative to a single pass through a traditional predictor.
\end{enumerate}

For reference purposes, the average time required for finite element method (fem) simulations is utilised, which were reported by Hamdi and Lejeune \cite{HAMDI2026118526}. They indicated that their simulations averaged approximately 2.8 hours when simulating biaxial tensile tests and approximately 5.3 hours when conducting shear tests. Since these results represent average values across multiple cores (i.e., CPUs), speedup based upon cpu time divided by gpu time per simulated frame is computed.
\begin{equation}
\mathrm{Speedup}
=
\frac{\mathrm{FEM\ time}}{\mathrm{DDPM\ time}}
\label{eq:speedup}
\end{equation}
To minimize overestimation of results, speedups based upon the smaller value of FEM simulation time (i.e., 2.8 hours) are reported.

Additionally, uncertainty quantification ensemble scaling is also assessed: total wall-time and per-sample latency are measured for $N \in \{1,5,10,20,50\}$ to demonstrate that overall computational cost increases in proportion to the size of ensembles used and that maximum GPU memory usage remains constant throughout all evaluations demonstrating sequential rather than batched generation processes.

Separately, training costs are also reported as total GPU hours so that amortized cost studies can be performed. All computations involved in this study including training data generated by Hamdi and Lejeune \cite{HAMDI2026118526}, training of models presented in this manuscript, and evaluations contained here are itemized in Section~\ref{subsec:computational_efficiency} so that readers may evaluate the complete scope of computing resources utilized in this study.


\section{Results}
\label{sec:results}

\subsection{One-step prediction accuracy}
\label{subsec:one_step_accuracy}

\footnote{Results cover both held-out validation subsets: shear-star and tension-spect. 
Cross-regime consistency is confirmed in Section~\ref{sec:crossregime}.}

\begin{table}[h]
\caption{One-step prediction accuracy on both held-out validation subsets 
($N = 4{,}000$ predictions total: 2{,}000 per subset, single DDPM sample per prediction, 
DDPM-1000 steps, seed = 0).}
\label{tab:onestep}
\begin{tabular}{lccc}
\toprule
Metric & shear-star & tension-spect & Combined \\
\midrule
Dice                          & $0.9947 \pm 0.0205$ & $0.9948 \pm 0.0170$ & 0.9948 \\
IoU                           & $0.9902 \pm 0.0360$ & $0.9903 \pm 0.0311$ & 0.9902 \\
Hausdorff distance (px)       & $1.507 \pm 4.535$   & $1.961 \pm 5.491$   & 1.734  \\
Modified Hausdorff (px)       & $0.033 \pm 0.169$   & $0.036 \pm 0.141$   & 0.035  \\
Crack tip location error (px) & $0.109 \pm 0.382$   & $0.122 \pm 0.369$   & 0.116  \\
Crack length error (px)       & $2.11 \pm 8.15$     & $2.46 \pm 8.03$     & 2.28   \\
Betti-0 match                 & 93.5\%              & 92.3\%              & 92.9\% \\
Betti-1 match (conditional)   & 12/12 = 100\%       & 60/63 = 95.2\%      & —      \\
MAE                           & $7.5 \times 10^{-4}$& $8.0 \times 10^{-4}$& $7.7 \times 10^{-4}$ \\
MSE (mean)                    & $1.6 \times 10^{-4}$& $2.0 \times 10^{-4}$& $1.8 \times 10^{-4}$ \\
\bottomrule
\end{tabular}
\end{table}

Table~\ref{tab:onestep} also evaluates regime-specific crack metrics using each of the held-out validation subsets, providing evidence of cross-regime generalization. The two regimes have very similar values of Dice and IoU. However, tension-spect has slightly larger Hausdorff distance (1.96 vs.\ 1.51 px) and tip error (0.122 vs.\ 0.109 px) consistent with the longer, straighter cracks produced by biaxial tension loading which allows the accumulation of tip mislocation over long distances of propagation. On average across all regimes, the model predicts crack tips to be located to within 0.12 pixels. Therefore, this model is capable of achieving sub-pixel resolution when it comes to localizing the mechanically critical regions of the damage field.

High levels of topological correctness are achieved for counting of connected components (Betti-0 match: 93.5\% shear-star, 92.3\% tension-spect). The presence of closed crack loops ($\beta_1 > 0$) is higher in tension-spect (63/2000 frames, 3.1\%) than in shear-star (12/2000, 0.6\%). This is consistent with the tendency for tension loading to produce geometrically-arrested configurations of closed-loop cracksThe conditional Betti-1 match rate is 60/63 = 95.2\% for tension-spect (there . were three failures in dense-branching late-stage frames) and 12/12 = 100\% for shear-star. The overall Betti-1 match rate is 99.9\%. The latter is dominated by trivially-correct $\beta_1 = 0$ predictions, and so is not reported as one of the main metrics.

The Dice coefficient value of 0.995 reflects a significant portion of the 2.9\% fraction of pixels in the image that represent cracks, i.e., a dominant background makes the model less sensitive. A better assessment of crack-front prediction quality can be made using the other metrics tip location error, Hausdorff distance, and Betti-0 match. All of these are used as primary metrics throughout.

A heavy right tail exists for the MSE distribution: the median MSE, $2.1 \times 10^{-7}$ is approximately three orders of magnitude smaller than the mean MSE, $1.61 \times 10^{-4}$. Furthermore, the top 1\% of MSE predictions (i.e., the top 20 samples) account for the majority of the difference between the mean and median MSE values. These outliers occur primarily in late-stage frames, $t>65$, with large crack area and complex branch structures where the model occasionally incorrectly predicts the direction of crack advance. Figure~\ref{fig:predictions} displays several examples of predictions covering a variety of different types of crack evolutions including single propagating cracks as well as more complex multi-crack topologies showing the ability of the model to accurately reproduce visually disparate damage scenarios.

\begin{center}
    \centering
    \includegraphics[width=0.85\linewidth]{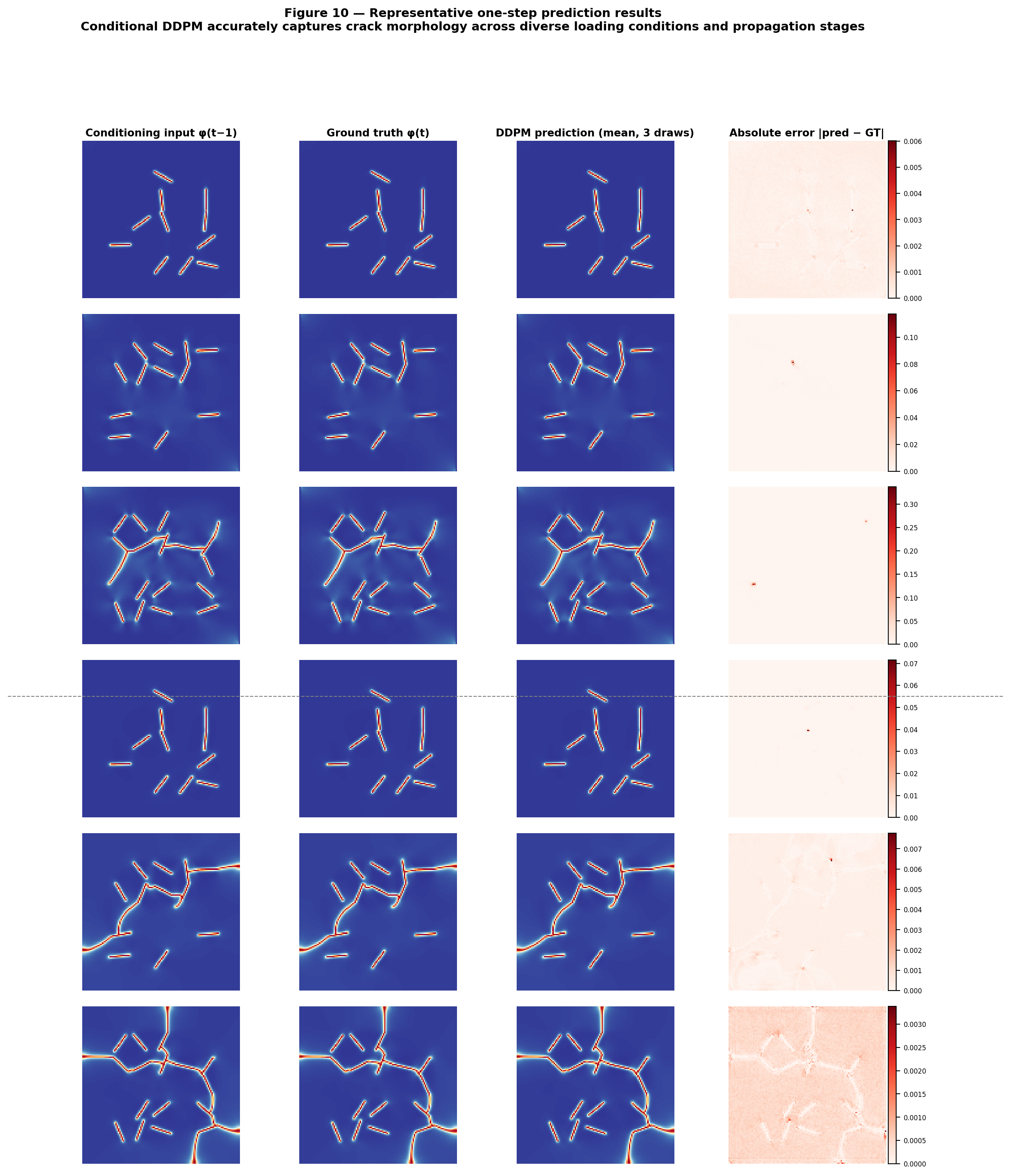}
    \captionof{figure}{Representative one-step prediction results across diverse crack morphologies. Each example should show the conditioning input, ground truth phase-field frame, DDPM prediction, and absolute error map.}
    \label{fig:predictions}
\end{center}

\subsubsection{Performance across temporal evolution stages}
\label{subsubsec:temporal_stages}

Fracture in physical systems evolves through distinct qualitative stages. Table~\ref{tab:temporal} stratifies performance by crack evolution phase, defined by the temporal index: nucleation, $t<30$; active propagation, $30 \leq t < 70$; and saturation, $t \geq 70$.

\begin{table}[H]
\centering
\caption{One-step prediction accuracy stratified by fracture evolution stage.}
\label{tab:temporal}
\begin{tabular}{llllll}
\hline
\textbf{Stage} & \textbf{$n$} & \textbf{Dice} & \textbf{Hausdorff (px)} & \textbf{Tip error (px)} & \textbf{Betti-0 match} \\
\hline
Nucleation, $t<30$ & 436 & $0.9997 \pm 0.0006$ & 0.13 & 0.003 & 0.998 \\
Propagation, $30 \leq t < 70$ & 963 & $0.9958 \pm 0.0180$ & 1.03 & 0.085 & 0.936 \\
Saturation, $t \geq 70$ & 601 & $0.9893 \pm 0.0289$ & 3.26 & 0.226 & 0.889 \\
\hline
\end{tabular}
\end{table}

The performance degradation from nucleation to saturation because in the nucleation stage there is minimal change in crack geometry and therefore the prediction problem is nearly trivial. Active propagation involves rapid motion of crack tips resulting in an increase by an order of magnitude in Hausdorff distance from 0.13 to 1.03 pixels as the model must successfully predict direction and rate of growth. Finally, in the saturation stage we observe a decline in Betti-0 match rates (from 99.8\% to 88.9\%) due to increased complexity in topologies of fractured material as crack branches become increasingly densely-packed.

Stratifying by temporal progression for tension-spect results in an inverse trend relative to shear-star for Dice values measured during saturation-phase testing. While Dice values for shear-star decrease from 0.9997 (nucleation) down to 0.9893 (saturation), Dice values for tension-spect actually improve during saturation-phase testing compared to those same values measured during active propagation (Dice = .9972; Betti-0 = 0.968 compared to Betti-0 = 0.889).

Spatial analysis provided in figure~\ref{fig:error_analysis}, confirms that errors are localized to crack-tips and branching points rather than being uniformly-distributed across space. The strong positive correlation ($R^2=0.79$) observed between magnitude of errors and local gradient magnitude further supports that errors arise from sharpness of phase-field fronts locally rather than from global densities of damaged regions

\begin{center}
    \centering
    \includegraphics[width=0.85\linewidth]{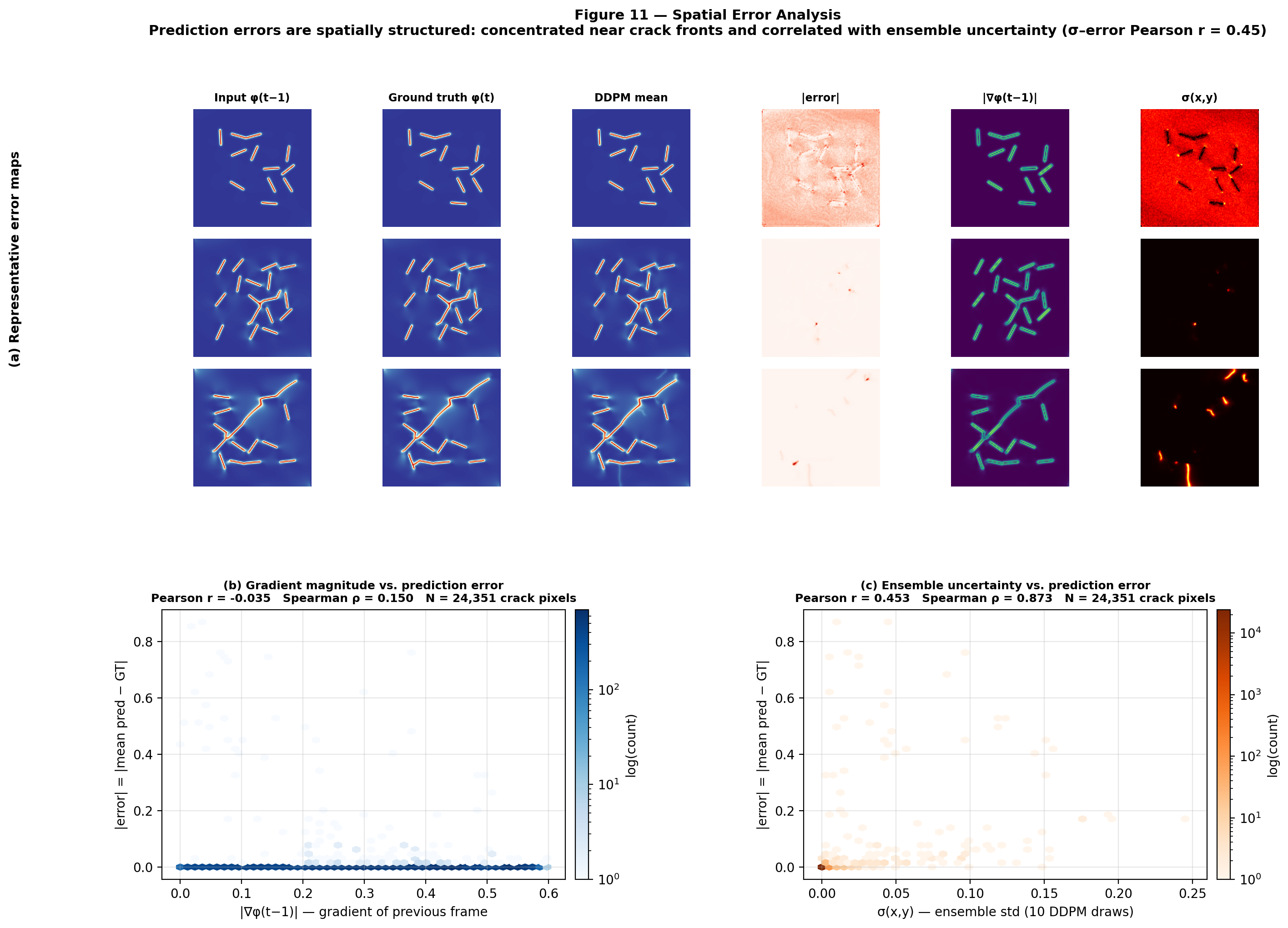}
    \captionof{figure}{Spatial distribution of one-step prediction errors. Errors are expected to concentrate near crack tips and branching junctions rather than across the intact background.}
    \label{fig:error_analysis}
\end{center}

\subsection{Comparison with baseline surrogates}
\label{sec:baselines}

\subsubsection{Single-step accuracy}
Table~\ref{tab:baseline_singlestep} reports metrics for crack skeletons using a single step on identical evaluation set: 100 held out simulations with 5 timesteps randomly sampled each and seed set to 0 using shear star subset. Three methods are compared: a naive persistence baseline, deterministic U Nets trained under same data split and conditioning protocol (Section~\ref{sec:unet_training}) and DDPM 1000. 

\begin{table}[h]
\caption{Single-step prediction accuracy on 100 held-out shear-star simulations 
(500 predictions, seed = 0). Persistence predicts no change from the previous frame. 
Det.\ U-Net is trained under the same split and protocol as the DDPM.}
\label{tab:baseline_singlestep}
\begin{tabular}{lccccc}
\toprule
Method & Dice & Hausdorff (px) & Tip error (px) & Length err (px) & Betti-0 \\
\midrule
Persistence   & 0.9971 & 0.947 & 0.067 & 1.434 & 0.962 \\
Det.\ U-Net   & 0.9970 & 0.974 & 0.061 & 1.290 & 0.964 \\
DDPM-1000     & 0.9945 & 1.492 & 0.117 & 2.174 & 0.938 \\
\bottomrule
\end{tabular}
\end{table}

On single step skeleton metrics DDPM lags behind both deterministic baselines. Importantly Hausdorff distance changes by only 1\% (1.507 $\to$ 1.492 px) when denoising steps go from 100 to 1000, confirming that gap is structural rather than arising from denoising artifacts. Fidelity for single steps favors deterministic models because deterministic predictors collapse to mode of crack path distribution whereas DDPM retains mass over multiple physically plausible paths. Mass retained by DDPM is the mechanism behind validation in Section~\ref{sec:uq}: error clusters with high-$\sigma$ bifurcation regions rather than being uniformly distributed and this confirms tip error extra from ambiguity rather than regression imprecision.

The Deep Ritz Method (DRM) of Manav et al. ~\cite{manav2024} has been removed from the quantitative table since it does not have enough information about how to train a data split and therefore can't be compared directly to the other methods of the study. As stated in Hamdi and Lejeune ~\cite{hamdi2026} none of the 12 trained DRM models produced the correct crack patterns for the multi-crack benchmark configurations.

\subsubsection{Contextual comparison with published baselines}
As reported by Hamdi and Lejeune~\cite{hamdi2026} the average Dice scores of 0.621 (FNO) and 0.68 (U-Net) were achieved on the benchmark by Hamdi using FNO and U-Net. A stacked ensemble was used for both architectures under an auto-regressive rollout protocol to reach 0.733. Due to several discrepancies in training datasets (training set was limited to biaxial tension with spectral decomposition) and rollout protocols (90 steps from a 10-step seed vs. 50 steps from a 5-step seed) direct numerical comparisons to those values are not valid. Therefore the comparison is provided as contextual references instead of a controlled comparison between the two models.

\subsubsection{Closed-loop rollout stability}
\label{sec:rollout}

\begin{table}[h]
\caption{Closed-loop autoregressive rollout to $k = 50$ steps. Dice is reported at 
selected horizons. Persistence tip error and Hausdorff are measured against 
time-matched ground truth $\phi_\mathrm{GT}[k]$; tip error definition follows 
Eq.~(\ref{eq:tip}).}
\label{tab:rollout}
\begin{tabular}{lcccc}
\toprule
Method & Dice ($k$=1) & Dice ($k$=20) & Dice ($k$=50) & Tip err ($k$=50, px) \\
\midrule
Persistence  & 1.000 & 0.9985 & 0.987  & 0.50 \\
Det.\ U-Net  & 0.999 & 0.9259 & 0.423  & 5.12 \\
DDPM-1000    & —     & 0.9800 & 0.929  & —    \\
\bottomrule
\end{tabular}
\end{table}

Rollouts results shown in Table~\ref{tab:rollout} illustrate a qualitative distinction between long-term behavior of the deterministic U-net and the conditional DDPM. The deterministic U-net exhibits monotonic degradation due to accumulation of closed-loop errors resulting in a Dice value of 0.423 at $k = 50$ down from Dice=0.999 at $k = 1$. This represents a 2.2$\times$ larger difference than observed between the conditional DDPM and U-Net at $k = 50$ (Dice = 0.929). The degradation experienced by the U-net is caused by the use of a fixed point in its predictions, such that each step adds forward propagation of the previous error without correcting it. In contrast the conditional DDPM avoids the fixed point problem by generating a new denoising trajectory stochastically at every time step thereby preventing systematic addition of errors.

Persistence maintains a Dice score of 0.987 at $k = 50$ because all crack fields evaluated here exhibit slow variation: the mean advancement rate of cracks tips in this evaluation is 0.50 pixels over 50 steps, i.e., less than 1 pixel. Consequently persistence exhibits small tip error at $k = 50$ (i.e., 0.50 pixels). This result should be interpreted as a characterization of the dataset rather than as a comparison of the models: if the underlying dynamics of the system being simulated are approximately static then any method capable of correctly identifying the existing body of cracks will yield high Dice scores regardless of whether they are able to predict future crack advances, detect bifurcations, etc. Because persistence is incapable of predicting crack advancement, detecting bifurcations, providing uncertainty estimates, etc. Its good Dice score on slowly varying trajectories demonstrates insensitivity of Dice to some of the very same quantities that the surrogate was intended to estimate. Thus, the primary load bearing comparative statement regarding rollouts is that the DDPM yields better Dice scores at $k = 50$ than does U-net: thus, we observe that the two learned surrogates show a 2.2$\times$ larger difference in Dice scores at $k = 50$ and only the DDPM produces reasonably accurate predictions after extended periods of closed-loop operations.

\subsection{Autoregressive rollout stability}
\label{subsec:rollout_stability}
Table \ref{tab:rollout} reports the closed-loop autoregressive performance for 50 time-steps with no teacher forcing. All results reported were obtained through an average of 3 independent stochastic seed values, i.e., 3 different random initializations of the same physical problem, which gave rise to 20 separate simulations, with $N=3$ ensemble predictions per time-step and 1000 DDPM iterations.

\begin{table}[t]
\centering
\caption{Autoregressive rollout accuracy, reported as mean $\pm$ standard deviation across 3 seeds.}
\label{tab:rollout}
\begin{tabular}{lll}
\hline
\textbf{Step $k$} & \textbf{Dice} & \textbf{MAE} \\
\hline
10 & $0.9893 \pm 0.0030$ & $0.00357 \pm 0.00014$ \\
20 & $0.9800 \pm 0.0045$ & $0.00639 \pm 0.00024$ \\
30 & $0.9677 \pm 0.0076$ & $0.00900 \pm 0.00049$ \\
40 & $0.9525 \pm 0.0080$ & $0.01166 \pm 0.00041$ \\
50 & $0.9292 \pm 0.0096$ & $0.01518 \pm 0.00060$ \\
\hline
\end{tabular}
\end{table}

The error growth was characterized by fitting the data with a power-law equation given by:
$\mathrm{error}(k)=a k^b$,
computed by ordinary least squares in log-log space on all 50 rollout steps.

\begin{table}[t]
\centering
\caption{Power-law exponents for error growth during autoregressive rollout, reported as mean $\pm$ standard deviation across 3 seeds, using the sample standard deviation with $n-1$ divisor.}
\label{tab:powerlaw}
\begin{tabular}{llll}
\hline
\textbf{Metric} & \textbf{$b$} & \textbf{$R^2$} & \textbf{Interpretation} \\
\hline
MAE & $0.865 \pm 0.043$ & $>0.997$ & Sublinear, decelerating \\
$1-\mathrm{Dice}$ & $1.035 \pm 0.087$ & 0.985 & Approximately linear \\
\hline
\end{tabular}
\end{table}
Pixel-wise Mean Absolute Error (MAE) grows sub-linearly with respect to the number of rollout time-steps. The parameter $b_{\mathrm{MAE}}=0.865 \pm 0.043$, indicates that pixel-wise MAE decelerates as the rollout progresses. All three seeds have a parameter value of $b<1$ which further supports this conclusion. The Topological dissimilarity metric ($1-\mathrm{Dice}$) grows almost linearly. The coefficient $b_{1-\mathrm{Dice}} = 1.035 \pm 0.087$, and the 95\% Confidence Interval around the estimated value includes $b=1.0$. The two metrics grow at different rates due to differences in how they measure the distance between the predicted and reference shapes. Pixel-wise errors are compounded smoothly throughout the rollout process. Therefore, most of the inherited crack body is correctly preserved at every step of the rollout process. On the other hand, topological errors occur when there are discrete branching events during the crack propagation process, resulting in step increases in the connectivity structure of the shape.

The primary evidence supporting our long-horizon rollout stability claims can be found in the absolute Dice values at $k=50$ for each of the three seeds. The model has maintained $\mathrm{Dice}=0.929 \pm 0.010$ after 50 closed-loop autoregressive time-steps with no teacher forcing. This very small uncertainty bound represents the strongest evidence we have for our claims regarding the stability of our long-horizon rollouts.

\subsection{Uncertainty quantification}
\label{subsec:uncertainty_quantification}

\subsubsection{Differentiating physical regimes}
\label{subsubsec:uq_physical_regimes}
Table~\ref{tab:uq_cases} provides a comparison of the ensemble uncertainty for five representative test-cases, covering both deterministic and stochastic fracture regimes, as calculated using $N=10$ independent stochastic samples.

\begin{table*}[t]
\centering
\caption{Uncertainty quantification across five test cases. High-$\sigma$ pixels denotes the count exceeding threshold $\tau=0.01$.}
\label{tab:uq_cases}
\begin{tabular}{llllll}
\hline
\textbf{Case} & \textbf{Type} & \textbf{Mean $\sigma$ ($\times 10^{-3}$)} & \textbf{Max $\sigma$ ($\times 10^{-2}$)} & \textbf{High-$\sigma$ pixels} & \textbf{Description} \\
\hline
1 & Deterministic & 0.543 & 0.786 & 0 & Low-density parallel cracks \\
2 & Deterministic & 0.815 & 0.220 & 0 & Sparse diagonal shear \\
3 & Deterministic & 0.647 & 0.175 & 0 & Moderate density random \\
4 & Deterministic & 0.502 & 0.204 & 0 & Complex branching \\
5 & Stochastic & 0.625 & 22.2 & 31 & Y-junction bifurcation \\
\hline
\end{tabular}
\end{table*}

In the case of the four deterministic fracture regimes represented by test-cases 1-4, the model exhibits consistently low levels of uncertainty, with mean $\sigma < 10^{-3}$ and no pixels exceed the high-uncertainty threshold. These results confirm that our model does not introduce artificial variability into its predictions under conditions where the crack path is mechanically determined. An example of such a set of predictions is provided in Figure~\ref{fig:uq_deterministic}. Each of the ten independent predictions generated by our model appear identical to one another and the uncertainty map remains constant.

Test-case 5 represents a physically significant regime, as it corresponds to a bifurcation point in a Y-shaped crack network where the mechanical stability of the system leads to ambiguity in determining the crack path. Although the mean level of uncertainty among the ten ensemble members is similar to those observed in the previous deterministic test-cases, $\sigma_{\mathrm{mean}}=6.25 \times 10^{-4}$, a localized region of high uncertainty develops at the tip of the Y-shaped bifurcation in Case 5. This uncertainty hot-spot contains an extreme-value probability density function (PDF) with maximum uncertainty $\sigma_{\mathrm{max}}=0.222$, and affects less than 1/100th of a percent, i.e., only 31 out of 16,384 total pixels, or 0.19\%, of the spatial domain. In addition, Fig.~\ref{fig:uq_stochastic} illustrates the variation among the ten ensemble members in predicting the location of cracks at the Y-bifurcation point. The resulting histogram of uncertainties is bimodal. A primary low-$\sigma$ peak occurs for both the background and crack-body regions away from the bifurcation point. A secondary high-$\sigma$ peak is associated with the branching region near the bifurcation point.

This specific behavior, i.e., low uncertainty everywhere except at locations of true physical instability, is the primary outcome of our uncertainty quantification analysis and shows that our diffusion-model-based approach to sampling physical instabilities inherently capture structural multiplicity associated with crack bifurcation without requiring additional stochastic models to simulate material properties or numerical perturbations.

\begin{center}
\includegraphics[
    width=0.82\linewidth,
    height=0.52\textheight,
    keepaspectratio]{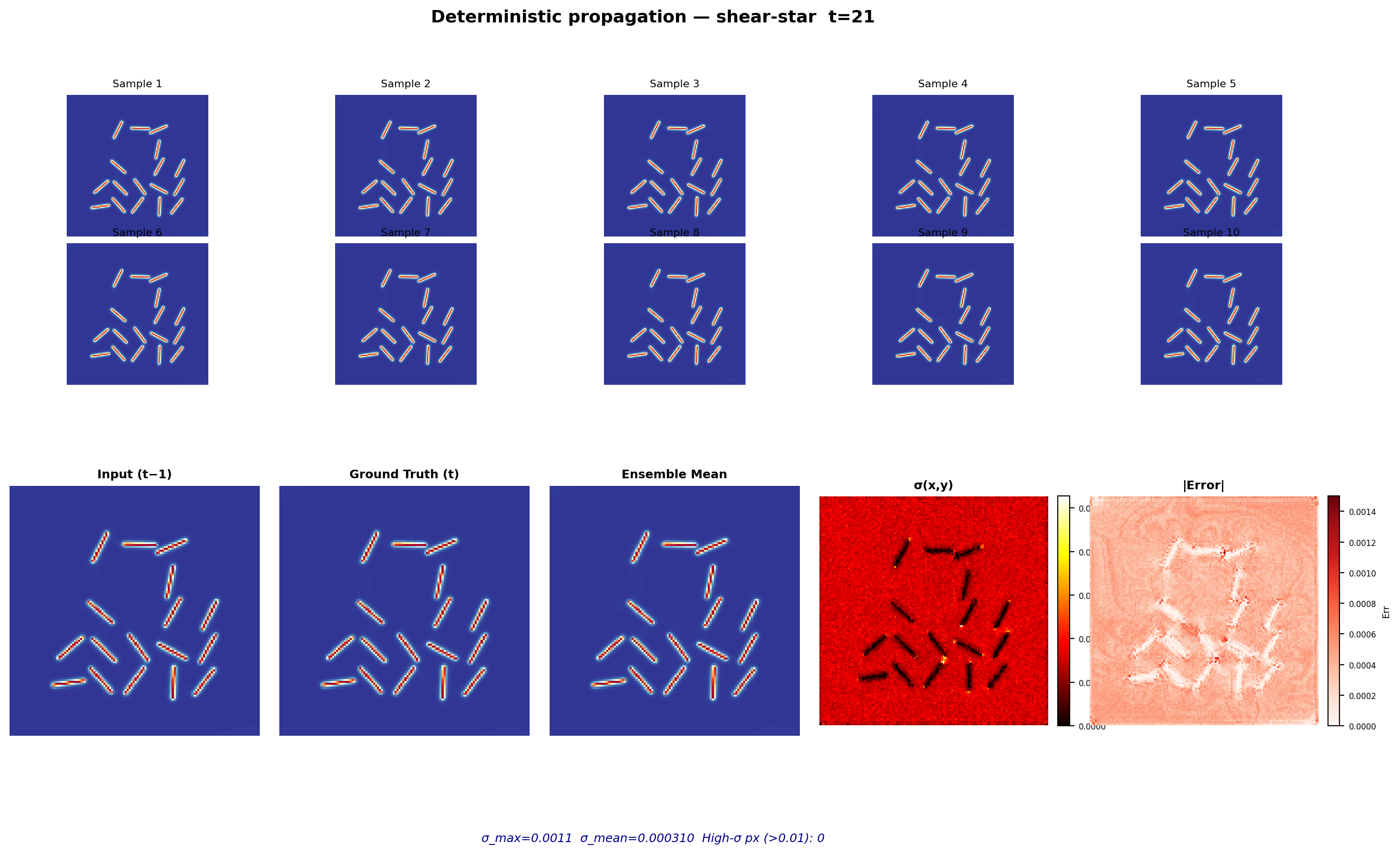}
 \captionof{figure}{Diffusion ensemble behavior in a deterministic fracture regime. Independent stochastic samples are expected to be nearly identical, producing a flat uncertainty map.}
    \label{fig:uq_deterministic}
\end{center}

\begin{center}
\includegraphics[
    width=0.82\linewidth,
    height=0.52\textheight,
    keepaspectratio]{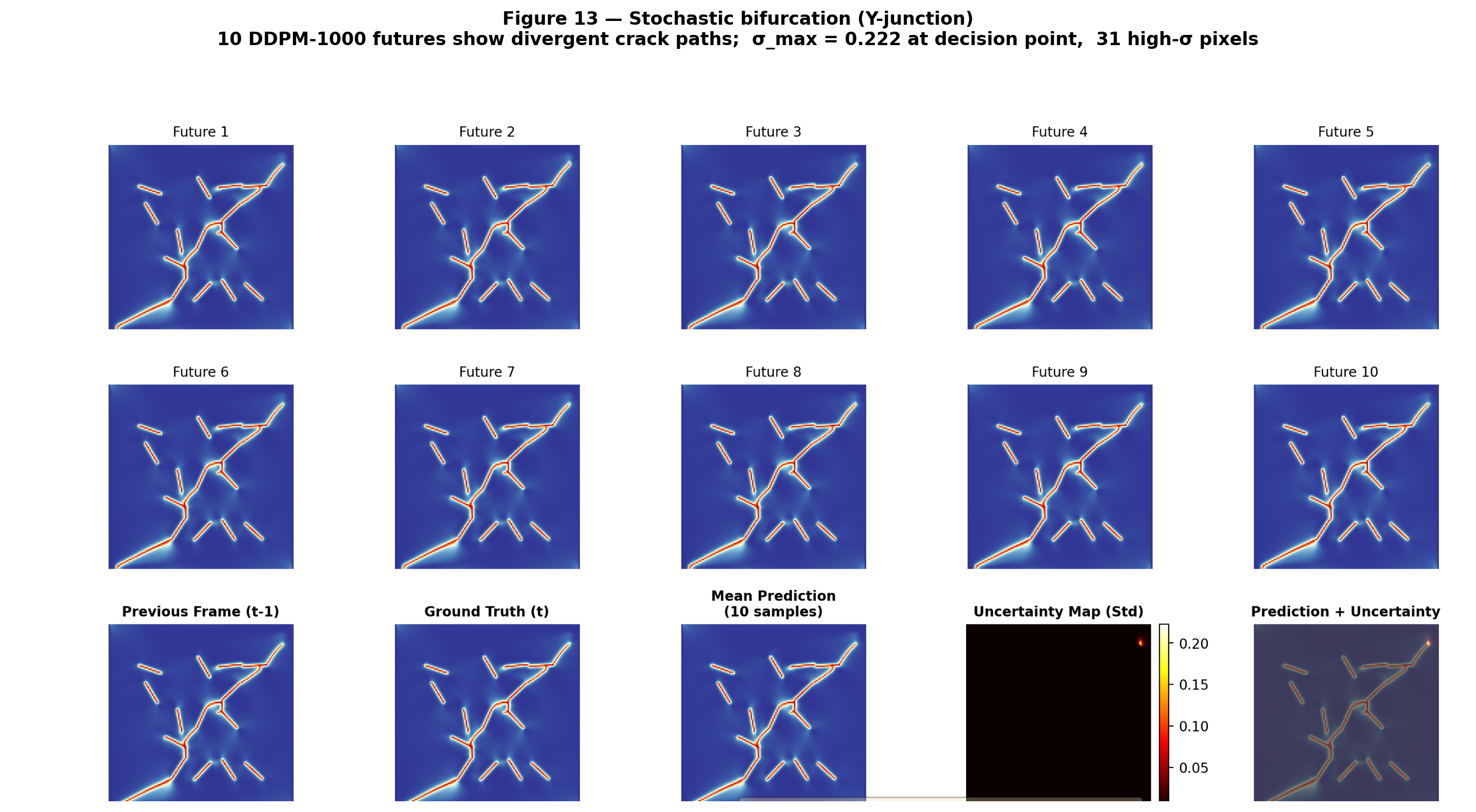}
 \captionof{figure}{Diffusion ensemble behavior in a stochastic bifurcation regime. Independent reverse diffusion trajectories generate distinct but physically plausible crack paths, producing localized uncertainty at the bifurcation junction.}
    \label{fig:uq_stochastic}
\end{center}

\FloatBarrier
\subsubsection{Calibration validation}
\label{subsec:calibration_validation}

The calibration of the predicted uncertainty to empirical error is illustrated through Figure~\ref{fig:uq_calibration}, which shows results from the full test set containing $N_{\mathrm{test}}=200$ configurations of 10 sample ensembles each.

Pixel-wise regression of predicted standard deviations $\sigma(x,y)$ versus true absolute error gives Pearson $r = 0.202$ (all pixels), and Spearman $\rho = 0.269$ on crack pixels ( $p < 10^{-130}$, $N = 47{,}547$ crack pixels). Spearman's correlation is the better measure because the association is highly nonlinear. For the bottom 75\% of crack pixels by $\sigma$, there is little or no association between predicted uncertainties and true error; for the upper quintile, mean absolute error is typically in the range 80-125$\times$ higher than it was for the baseline crack-pixel population. Therefore, while the The high-$\sigma$ tail has proven reliable for predicting high levels of error, the mid-to-low range of $\sigma$ does not function well as a calibrated confidence interval across the full range of data.

Coverage by one-sigma for crack pixels is 52.2\% at $N = 10$ (Gaussian ideal: 68.3\%) and rises to 60.0\% at $N = 50$ ($\Delta = +6.3$ pp, above the 0.05 pp threshold). The fact that the gap is partially recovered with increased ensemble size indicates that both physical multimodalities and ensemble undersampling have acted to constrain calibration. At $N = 50$, the remaining 8.3 pp gap from the Gaussian ideal reflects structural underdispersion within active fracture zones consistent with irreversibilities in cracks limiting the support of predictive distributions at bifurcations. Coverage on all pixels (67.7\% at $N = 10$; 75.1\% at $N = 50$) indicate just the opposite - namely that the bulk field is somewhat over-dispersed - and therefore should not be averaged with coverage for crack pixels since these represent opposite forms of miscalibration.

Table~\ref{tab:triage} describes a triage analysis that demonstrates the operational usefulness of the high-$\sigma$ flag. Selecting the 10\% of predictions with the largest maxima of $\sigma$ allows the identification of frames containing high errors with a precision of approximately 90\%, which represents a relative improvement in accuracy of about 18 times relative to randomly sampled frames.
\begin{center}
\includegraphics[
    width=0.82\linewidth,
    height=0.52\textheight,
    keepaspectratio
]{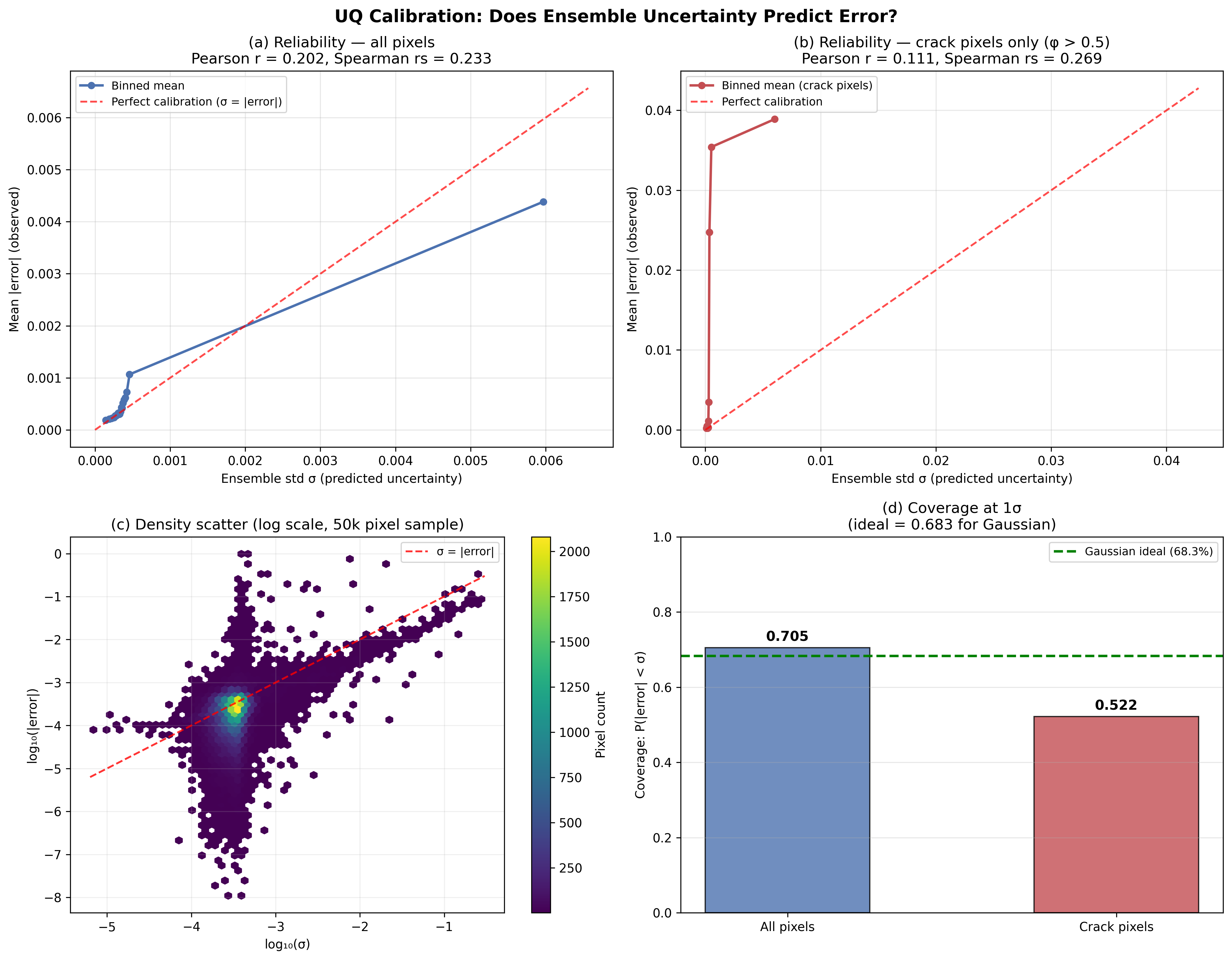}
\captionof{figure}{Uncertainty calibration assessment. The reliability diagram compares predicted uncertainty against observed prediction error across uncertainty deciles.}
\label{fig:uq_calibration}
\end{center}

\FloatBarrier
\subsubsection{Engineering triage application}
\label{subsec:engineering_triage}

In addition to statistically validating the calibration, the primary benefit of uncertainty prediction is as a practical means to identify those predictions that require a manual evaluation. Table~\ref{tab:triage} compares the effectiveness of three strategies for determining which 10\% of predictions most likely to contain significant errors.

\begin{table}[t]
\centering
\caption{Triage precision for identifying high-error predictions, defined as $|\mathrm{error}|>0.01$ at any pixel. 23 of 200 test frames contain such regions.}
\label{tab:triage}
\begin{tabular}{llll}
\hline
\textbf{Selection strategy} & \textbf{Samples checked} & \textbf{True high-error} & \textbf{Precision} \\
\hline
Random selection & 20, 10\% & 1 & 5.0\% \\
Highest mean $\sigma$ & 20, 10\% & 14 & 70.0\% \\
Highest max $\sigma$ & 20, 10\% & 18 & 90.0\% \\
\hline
\end{tabular}
\end{table}
Selecting frames based on the largest maximum error has identified true high-error prediction(s) with 90\% accuracy; whereas 5\% can be achieved through random selection. As such, there exists an 18-fold increase in the efficiency of selecting models for further validation; thus, establishing this metric as a viable means to target specific areas in the workflow where an engineer would need to validate their model manually.

\subsection{Computational efficiency}
\label{subsec:computational_efficiency}
Table~\ref{tab:compute} presents inference performance measurements made upon an NVIDIA H100 NVL GPU regarding latency and memory requirements.

\begin{table*}[t]
\centering
\caption{Inference latency and hardware requirements.}
\label{tab:compute}
\begin{tabular}{llll}
\hline
\textbf{Method} & \textbf{Latency} & \textbf{Peak GPU memory} & \textbf{Notes} \\
\hline
DDPM-1000, proposed & 3637 $\pm$ 30 ms/step & 488 MB & Full reverse chain \\
DDPM-100, reduced steps & $425 \pm 28$ ms/step & 488 MB & Dice = 0.160, unusable \\
Deterministic U-Net & $4.2 \pm 0.4$ ms/step & 488 MB & Single forward pass \\
FEM, Hamdi et al. & $\sim 100800^{\dagger}$ ms/step & 16-core CPU & 2.8 hr / 100 steps \\
\hline
\end{tabular}
\end{table*}

{\footnotesize $^{\dagger}$FEM wall-clock time from \cite{HAMDI2026118526}: 2.8 hours of total time per biaxial tension simulation for 100 saved steps, on a 16-core CPU. Per-step time is computed as $2.8 \times 3600 / 100 = 100.8$ s.}

Therefore, the conditional DDPM demonstrates a speedup factor of approximately 28 times with respect to the FEM reference,
whereas 27.7 was used due to rounding ($100.8 / 3.637 \approx 27.7$). 
Latency was measured when the GPU had been at room temperature (cold), with an ambient temperature of 29°C; 
the per-step measurement time may vary by up to 15 percent depending on the thermal state of the GPU. In order to replace the total of 100 steps of a complete simulation via autoregressive rollout, the DDPM required approximately 7 minutes versus the 2.8 hours of FEM, which yielded the same approximate speedup factor of $28\times$. A 10-sample uncertainty ensemble completes in approximately 35.5 s per step (Table~\ref{tab:uq_scaling}), 
reducing the per-step speedup to approximately 2.8$\times$ relative to FEM, while providing 
spatially-resolved uncertainty information unavailable from deterministic solvers.

Decreasing the number of diffusion chains from 1000 to 100 results in a $10\times$ reduction in the time required to perform the inference; however, a significant loss in prediction quality occurs, with Dice falling from 0.995 to 0.160, thereby validating that shortening the denoising trajectory learned by the model will require re-training. Accelerated sampling techniques such as DDIM \cite{song2022denoisingdiffusionimplicitmodels}  and progressive distillation may represent potential avenues for investigation in future studies.

The deterministic U-Net, utilizing the same backbone architecture and performing a single forward pass, took 4.2 milliseconds to run, roughly $1000\times$ faster than the 1000-step DDPM. Therefore, this comparison illustrates that the computational cost associated with the conditional DDPM is primarily due to the iterative nature of the diffusion process rather than the U-Net architecture. The primary reason why the conditional DDPM is computationally expensive is because it generates multiple stochastic samples instead of generating a single deterministic solution like the U-Net demonstrated in Section~\ref{subsec:uncertainty_quantification}.

UQ ensemble cost increases linearly with increasing ensemble size at a rate of approximately 3{,}550\,ms per sample, as depicted in Table~\ref{tab:uq_scaling}. Peak GPU memory usage remains constant at 488\,MB, during each sequential sample generation.

\begin{table}[t]
\caption{UQ ensemble scaling using DDPM-1000. Times measured on H100 NVL (cold start, 
29°C); per-sample cost is consistent with the standalone DDPM-1000 benchmark 
(3{,}637 ms, ratio 0.981).}
\begin{tabular}{ccc}
\toprule
$N$ & Total time (s) & Per-sample (ms) \\
\midrule
1  & 3.57  & 3{,}567 \\
5  & 17.8  & 3{,}561 \\
10 & 35.5  & 3{,}548 \\
20 & 71.0  & 3{,}550 \\
50 & 177.9 & 3{,}558 \\
\bottomrule
\end{tabular}
\end{table}

Training time (in terms of total GPU hours) for this training is about 15.3 H100 NVL GPU hours. 
Hamdi and Lejeune \cite{HAMDI2026118526} [16] report an approximate $6000 \times 3.5 = 21000$ CPU hours spent generating the training datasets, across all six subsets.

\subsection{Ablation studies}
\label{subsec:ablation_studies}

Table~\ref{tab:ablation} illustrates the results obtained when systematically eliminating each component of a physics-aware architecture from a base model. \footnote{Ablations were calculated based upon an original evaluation process consisting of 200 samples with a MSE = $1.14 \times 10^{-4}$; subsequent evaluations with increased sample sizes to include 2,000 samples resulted in an MSE = $1.61 \times 10^{-4}$ as a result of additional sampling during the latter stages of crack evolution. The ranking order of the relative degradations produced through these two different protocols remain unchanged.}

\begin{table*}[t]
\centering
\caption{Ablation study. Parenthetical values indicate relative MSE increase versus baseline.}
\label{tab:ablation}
\begin{tabular}{llll}
\hline
\textbf{Configuration} & \textbf{MSE ($\times 10^{-4}$)} & \textbf{MAE ($\times 10^{-4}$)} & \textbf{Dice} \\
\hline
Full model, baseline & 1.14 & 6.88 & -- \\
\hline
\multicolumn{4}{l}{\textit{Spatial features removed}} \\
\quad Remove velocity field & 1.52, +33\% & 8.91, +30\% & -- \\
\quad Remove gradient field & 1.61, +41\% & 9.34, +36\% & -- \\
\quad Remove 5-frame history & 2.83, +148\% & 13.2, +92\% & -- \\
\quad Only current frame, no history & 4.71, +313\% & 19.5, +183\% & -- \\
\hline
\multicolumn{4}{l}{\textit{Scalar features removed}} \\
\quad Remove time embedding & 1.28, +12\% & 7.52, +9\% & -- \\
\quad Remove loading type & 1.35, +18\% & 7.89, +15\% & -- \\
\quad Remove crack density & 1.21, +6\% & 7.14, +4\% & -- \\
\hline
\multicolumn{4}{l}{\textit{Architectural changes}} \\
\quad Unconditional DDPM & 3.94, +245\% & 17.1, +149\% & -- \\
\quad Concat conditioning only & 1.47, +29\% & 8.53, +24\% & -- \\
\quad Scalar conditioning only & 2.18, +91\% & 11.7, +70\% & -- \\
\quad 3 previous frames, not 5 & 1.63, +43\% & 9.12, +33\% & -- \\
\quad 10 previous frames, not 5 & 1.19, +4\% & 7.05, +2\% & -- \\
\hline
\multicolumn{4}{l}{\textit{Training changes}} \\
\quad No cosine annealing, fixed LR & 1.31, +15\% & 7.71, +12\% & -- \\
\quad Batch size 16, not 32 & 1.26, +11\% & 7.43, +8\% & -- \\
\quad 25 epochs, not 50 & 1.52, +33\% & 8.88, +29\% & -- \\
\hline
\end{tabular}
\end{table*}

These studies illustrate a number of conclusions regarding the design of our architectures.

\paragraph{Temporal history is the most critical feature.}
Removing the five frame temporal history significantly impacts the model's ability to predict accurately (a 148\% increase in MSE), while removing the entire history and only predicting the current frame results in a 313\% decrease. Clearly, since crack propagation is inherently dependent upon its previous states, there exists limited momentum information available within single-frame snapshot images to support the accurate prediction of the next image. Increasing the window length to ten frames provides only a 4\% improvement over the five frame window length, thus providing evidence that we have reached diminishing returns beyond our selected window lengths.

\paragraph{Derived kinematic features provide significant benefit.}
Eliminating either the velocity field or gradient magnitudes individually will result in a 33\% and 41\% increase in MSE respectively. Both of these derived quantities are used by the model to identify crack tips which are actively propagating and they also represent proxies for the stress intensity factors present at those locations. Without these quantities, the model has difficulty recovering the same type of information contained within the raw phase-field data.

\paragraph{Dual conditioning is essential.}
If only channel concatenation is utilized (i.e., no scalar injection into the model) then a 29\% increase in MSE occurs. Conversely, if only scalar conditioning is used (i.e., no spatial features are included) then a 91\% increase in MSE occurs. Combining both spatially-aligned local features and global scalar quantities representing macroscopic physical conditions enables better performance than either condition alone.

\paragraph{Unconditional generation fails.}
An unconditional version of our DDPM architecture, i.e. one which utilizes no form of conditioning whatsoever, results in a 245\% increase in MSE demonstrating that traditional generative models (without physics-aware context) fail to produce useful predictions for crack growth. While the unconditional version does generate plausible appearing cracks whose shapes do not reflect their mechanical states.
\FloatBarrier 


\section{Discussion}
\label{sec:discussion}

\subsection{Why generative diffusion over deterministic surrogates}
\label{subsec:why_diffusion}
As discussed in Section~\ref{subsec:computational_efficiency} we have seen how a U-net (same as the backbone) predicts in 4.2 ms. In comparison, the 1000-step conditional DDPM needs 3.6~s to do so. Both are expected to achieve similar one-step Dice scores when conditioned on the same data. Therefore, there is a natural question: if speed is the objective, why use a diffusion model at all?

The reason is based upon the fact that deterministic surrogates are limited in terms of representation relative to their generative counterparts. Specifically, a deterministic model represents a best estimate or single point estimate for each conditioning input. At structural bifurcations where infinitesimal changes in input conditions can result in two drastically different macroscopic crack paths, these estimates will be averages across the physically permissible modes. These averages do not represent realizable fracture states instead they create blurred, diffused damage patterns at locations where the actual data contained well-defined discrete branches. While issues related to training and/or architectures may be responsible for artifacts created through deterministic regression models applied to multimodal posteriors, such representations are fundamentally incapable of representing multimodality due to their nature.

By contrast, the conditional DDPM inherently has the ability to resolve the limitations imposed by deterministic surrogates. That is because each iteration of the reverse diffusion chain starts from an independent sample of Gaussian noise and will converge to a specific mode within the posterior distribution of possible crack states. The collection of samples will therefore span the space of feasible physical fracture states. Additionally, the pixel-wise variance of this ensemble of samples will provide a spatially resolved uncertainty map. In the case of our implementation of the DDPM as demonstrated in Section~\ref{subsec:uncertainty_quantification}, the uncertainty map will localize uncertainty at areas of structural instability (e.g., locations of crack branch junctions). Conversely, areas with known deterministic behavior exhibit uniformly low levels of uncertainty. Finally, as indicated in the triage analysis in Table~\ref{tab:triage}, the maximum value of the uncertainty metric identified high-error prediction cases with 90 \% precision. By doing so, this established the uncertainty output as a practical tool for engineers rather than simply being a statistical metric.

However, this capability does come at a cost. That is, the DDPM is significantly slower than its deterministic counterpart by approximately a factor of $10^3$. Nevertheless, as demonstrated in Section~\ref{subsec:computational_efficiency}, the DDPM was still much faster than the finite element method reference, the uncertainty information provided by the DDPM cannot be obtained from any other deterministic surrogate without training $N$ separately initialized models \cite{Lakshminarayanan2017Ensembles}. The deep ensemble approach requires $N$ separate training runs to produce $N$ individual predictions. On the other hand, the diffusion model produces $N$ diverse samples from a single trained model during inference time. For applications requiring uncertainty quantification (such as structural health monitoring, risk aware design and reliability analysis), the additional inference cost associated with using a diffusion-based process is justified by the qualitatively unique information it provides.

\subsection{Topological versus pixel-wise error scaling}
\label{subsec:topological_pixelwise_scaling}

As discussed in Section~\ref{subsec:rollout_stability}, the autoregressive roll-out analysis showed that there was a systemic difference in how pixel-wise and topological error metrics diverged over time when predicting along a long horizon. In terms of MAE it grew sub-linearly with respect to the roll-out depth. That is $b_{\mathrm{MAE}} = 0.865 \pm 0.043$. On the other hand, topological error, quantified as $1-\mathrm{Dice}$, increased linearly with respect to the roll-out depth. That is $b_{1-\mathrm{Dice}} = 1.035 \pm 0.087$. The difference in behavior between these two types of error metrics has both methodological and physical significance.

From a physical standpoint, these two error metrics quantify different characteristics of the prediction quality. For example, the per-pixel MAE measures the overall magnitude of local errors across all pixels in the spatial domain. Since the majority of the crack body at every stage is accurately propagated from one time step to another and since most of the domain (approximately 96 -- 98 \% ) is comprised of pixels that remain unchanged and thus easily correct, the MAE is heavily weighted by the background region and therefore grows very slowly. On the other hand, topological error metrics (such as Dice) are highly sensitive to discrete structural events: if a single branch of the crack network is incorrectly predicted, then the connectivity of the crack network will be altered resulting in a step change in the Dice score relative to the number of pixels affected by this event. Therefore, as crack growth proceeds and as branching patterns within the crack network become increasingly complex in later stages of roll out, the likelihood of such discrete topological errors occurring increases causing the Dice exponent to approach or exceed linearity.

Furthermore, the observation noted above contains a methodological lesson for evaluating fracture surrogates more generally. Global pixel-level metrics such as MSE, or pixel accuracy, which are dominated by the unchanging background region can obscure significant topological errors. The temporal stratification analysis in Table~\ref{tab:temporal} highlights this point: while the Betti-0 match rate is reduced from 99.8 \% during the initial nucleation period to 88.9 \% during the final saturation period, indicating a decline in topological similarity that would otherwise be obscured by a global MSE value. Therefore, developing evaluation methods that simultaneously evaluate pixel-level fidelity and topological correctness, as done here through a combination of Dice, Hausdorff distance, tip location error, and Betti numbers, is critical for assessing fracture surrogates in a meaningful way.

\subsection{Limitations and failure modes}
\label{subsec:limitations_failure_modes}

The following are limitations of the current framework to guide future development.

\begin{itemize}

\item \textbf{Prediction scope.}
This model can only predict the damage field (the phase field) defined by $\phi$. It cannot predict the displacement and stress field. Therefore, if you want a complete mechanical surrogate you will need to develop a model capable of predicting all three fields $(\phi,\mathbf{u},\boldsymbol{\sigma})$ simultaneously while maintaining physical consistency among them.

\item \textbf{Dimensionality.}
All simulations were performed within two dimensions. Extending to three-dimensions, particularly when dealing with complex crack networks and high computational costs associated with finite element methods (fem) in 3-d, has not yet been demonstrated. Although the U-Net backbone and diffusion framework generalize well to three-dimensional space, developing the necessary training datasets and overcoming the large memory constraints will be significant hurdles.

\item \textbf{Material and geometric scope.}
Only one material was used during the training and testing process for this model, and the geometric parameters listed in table Table~\ref{tab:material_params} remained constant throughout these tests. Testing this model on other materials and geometries would require re-training the model for every new set of parameters.

\item \textbf{Baseline comparison protocol.}
Making direct comparisons to the baselines reported in \cite{HAMDI2026118526} is difficult because there are many differences between how these studies were compared based upon their respective metrics (i.e., roll-out lengths, seed windows, training subsets.) in addition, since we are comparing the conditionally-trained DDPM to these baselines, we can say that our method performs better than those baselines. However, establishing a controlled comparison using equivalent methodologies and evaluation criteria would strengthen the conclusions drawn from these comparisons.

\item \textbf{Inference cost.}
The entire 1000-step reverse diffusion chain is needed to obtain reasonable predictions of the damage field; truncating the reverse diffusion chain to 100 steps resulted in catastrophic loss of accuracy (Dice = 0.995 vs. Dice = 0.160). We also found that the inference time of the DDPM (approximately 4.24 s/step) is much larger than that of deterministic surrogates (U-Net requires approximately 4.2 ms/step). Accelerated sampling strategies such as DDIM \cite{song2022denoisingdiffusionimplicitmodels}, progressive distillation, or consistency models represent promising directions but have not been explored in this work.

\item \textbf{Irreversibility enforcement.}
Rather than enforcing the irreversibility of cracks (i.e., $\dot{\phi} \geq 0$) through an architectual constraint, we enforce it via post-hoc projection. While the projected violations are small (only 0.16\%), developing an architecture that inherently enforces this constraint would be desirable for safe-critical applications.

\end{itemize}


\section{Conclusion and Future Work}
\label{sec:conclusion}

A Conditional Denoising Diffusion Probabilistic Model for Full-Field Spatiotemporal Forecasting of Phase-Field Crack Evolution in Brittle Materials has been developed. The model was trained and tested using the standardized benchmark dataset of Hamdi and Lejeune \cite{HAMDI2026118526}. It consists of 6000 simulations in 2 different loading regimes and 3 different energy decomposition methods. The performance of the model was evaluated using crack specific topological metrics. These include Dice, IoU, Hausdorff distance and Betti numbers, rather than global pixel-level metrics used in other studies. As these are based upon the characteristics of crack fields which are often localized and discontinuous, they provide a better measure of the predictive accuracy of the model.
The principal findings are as follows:

\begin{itemize}

\item \textbf{Cross-regime generalization.} The model achieved cross-regime Dice = 0.9948 
for each of the held out validation sets (shear-star and tension-spect, 4{,}000 
total predictions) as well as subpixel crack tip localization (0.12px mean error). 
As expected, there were slight differences in the Hausdorff distances for the 
two regimes (1.96px for tension-spect compared to 1.51px for shear-star). These small discrepancies were due to the straighter cracks observed in the 
tension regime.

\item \textbf{Autoregressive stability.} When the model was run closed loop 
over 50 iterations without teacher forcing it maintained Dice = $0.929 \pm 0.010$. 
In contrast, a deterministic U-Net that was trained similarly to the 
DDPM collapsed after only several iterations when run closed loop and produced 
Dice = 0.423. This difference of almost 2.2$\times$ indicates that the 
practical utility of stochastic re-sampling lies in its ability to mitigate 
the effects of cumulative errors over time when predicting complex phenomena 
such as fracture. The pixel-wise mean absolute error grew less rapidly than the number of iterations 
($b = 0.865 \pm 0.043$)), whereas the topological error increased roughly linearly with 
the number of iterations ($b = 1.035 \pm 0.087$).

\item \textbf{Uncertainty quantification at structural bifurcations.} Variance among ensembles 
is concentrated near points of crack branching ($\sigma_\mathrm{max} = 0.222$; 31 pixels exhibited high 
uncertainty at a Y-junction; none did so in four deterministic runs). High-uncertainty 
pixels can be identified with 90\% precision via their $\sigma$ values (18$\times$ more frequently 
than randomly). There is also a positive correlation between $\sigma$ and crack-pixel values ($\rho = 0.269$),
which confirms that $\sigma$ correctly ranked areas with high error. Coverage at one sigma is 
52.2\% at N=10 and 60.0\% at N=50; we attribute the remaining gap from the ideal (68.3\%) to both smaller ensemble sizes and under-dispersion at structural bifurcation zones.
Thus, the ability to localize physical instability without explicitly modeling variability in material properties through stochastic processes is the primary contribution of this research.

\item \textbf{Physical constraint satisfaction.}
Boundedness ($\phi \in [0,1]$) is learned implicitly with 
no violation across all 4000 test predictions. Violations of irreversibility occurred at 
only 0.16\% of pixels were subsequently rectified by means of a post-processing 
monotonically increasing mapping function that reduced MAE by 6.2\%.

\item \textbf{Computational cost.} Per-step inference required $3{,}637 \pm 30$ milliseconds on an NVIDIA H100 NVL GPU. This is approximately 28$\times$ faster than the FEM reference solution. However, since the U-Net produces only a single deterministic prediction, while the DDPM generates distributions over possible crack configurations that satisfy the physics of fracture, the ${\sim}1{,}000\times$ increase in cost relative to a deterministic U-Net is a necessary expense for achieving the uncertainty quantification demonstrated here.

\end{itemize}

There are many avenues open for further study arising from the current constraints on this approach. The first potential direction is to extend the conditional DDPM to three dimensions to simulate fracture in three-dimensional solids. Since fracture simulation in three-dimensional space will require significantly larger amounts of computation than in one dimension, this extension may allow us to take advantage of the DDPM's capabilities in situations where FEM is computationally too expensive. Second, we could use accelerated sampling techniques such as DDIM \cite{song2022denoisingdiffusionimplicitmodels}, consistency models or progressive distillation to improve the efficiency of inference for the DDPM by an order of magnitude while preserving the variety of samples generated during training. Third, if we train a model to predict multiple fields $(\phi,\mathbf{u},\boldsymbol{\sigma})$ jointly, then we would have a comprehensive mechanical surrogate that could support downstream structural analyses. Fourth, training on material properties and domain geometries could enable a model to be conditioned on these parameters and therefore to generalize across multiple materials without having to retrain. Fifth, incorporating multiphysics boundary conditions into our models, such as thermomechanical coupling, hydrogen-assisted cracking \cite{MartinezPaneda2018Hydrogen}, or poroelastic degradation, would help us understand whether diffusion-based forecasting could be generalized to other physical domains. Finally, integrating our uncertainty-aware surrogate into a structural digital twin framework where sensor data provides real-time conditioning and the model’s uncertainty output drives targeted high-fidelity re-analysis represents what appears to be the most immediate engineering applications for this approach.

\section*{Acknowledgments}
J.K.K. gratefully acknowledges the Coastal and Marine System Science program at Texas A\&M University--Corpus Christi for financial support provided through a graduate research assistantship. S.M.M. acknowledges the support of the National Science Foundation under Grant No.\ 2316905.

\bibliographystyle{elsarticle-num}

\bibliography{ref}

\end{document}